\documentstyle[11pt,bezier]{article}
\input{amssym.def}
\input{amssym}
\newtheorem{th}{\arabic{section}.\arabic{abz}. Theorem}

\newtheorem{defi}{\arabic{section}.\arabic{abz}. Definition}

\newtheorem{prop}{\arabic{section}.\arabic{abz}. Proposition}

\newtheorem{cor}{\arabic{section}.\arabic{abz}. Corollary}

\newtheorem{conv}{\arabic{section}.\arabic{abz}. Convention}

\newtheorem{exa}{\arabic{section}.\arabic{abz}. Example}

\newtheorem{abza}{\arabic{section}.\arabic{abz}.}

\newtheorem{rem}{\arabic{section}.\arabic{abz}. Remark}

\newcounter{abz}[section]
\newcounter{sabz}[abz]
\addtocounter{section}{-1}
\newcommand{\abz}{\refstepcounter{abz}}
\newcommand{\sabz}{\refstepcounter{sabz}}

\newcommand{\con}{\rule{5pt}{1pt}\rule{1pt}{5pt}\,}
\newcommand{\lbr}{\linebreak[0]}
\def\d{\partial}
\def\Ci{C^{\infty}}
\def\G{\Gamma}
\def\la{\lambda}
\def\La{\Lambda}
\def\ra{\mathop{\rm rank}\nolimits}
\def\cora{\mathop{\rm corank}\nolimits}
\def\ke{\mathop{\rm ker}\nolimits}
\def\Re{\mathop{\rm Re}\nolimits}
\def\Im{\mathop{\rm Im}\nolimits}
\def\Lie{\mathop{\rm Lie}\nolimits}
\def\Vect{\mathop{\rm Vect}\nolimits}
\def\Span{\mathop{\rm Span}\nolimits}
\def\Sing{\mathop{\rm Sing}\nolimits}
\def\im{\mathop{\rm im}\nolimits}
\def\i{{\rm i}}
\def\codim{\mathop{\rm codim}\nolimits}
\def\diam{\mathop{\rm diam}\nolimits}

\def\E{{\cal E}}
\def\R{{\Bbb R}}
\def\C{{\Bbb C}}
\def\A{{\frak g}}
\def\D{{\frak g}^*}
\def\O{{\cal O}}
\def\K{{\cal K}}
\def\J{{\cal J}}

\begin{document}
\title{Symplectic realizations of bihamiltonian structures}
\author{Andriy Panasyuk\thanks{Partially supported
by the Polish grant KBN 2 PO3A 135 16.} \\ Division of
Mathematical Methods in Physics, \\ University of Warsaw,
\\ Ho\.{z}a St. 74, 00-682 Warsaw, Poland, e-mail:
panas@fuw.edu.pl\\ and \\ Mathematical Institute of Polish
Academy of Sciences\\ \'{S}niadeckich St. 8, 00-950 Warsaw,
Poland} \date{}

\maketitle

\hspace{6.5cm}{\em To my wife Larysa}

\section{Introduction}

A smooth manifold $M$ is endowed by a Poisson
pair if two linearly independent bivector fields
$c_1,c_2$ are defined on $M$
and moreover  $c_\la
=\la_1c_1+\la_2c_2$ is a  Poisson tensor for any $\la
=(\la_1,\la_2)\in\R^2$. A bihamiltonian structure
$J=\{c_\la\}$  is the whole 2-dimensional family of
tensors.

There are two classes of bihamiltonian structures playing
important role in the theory of completely integrable
systems. The geometries corresponding to these classes are
different and so are the ways for appearing of functions in
involution.

The first class, called symplectic in this paper, is
characterized by the condition that $\ra c_\la=\dim M$
(cf. Convention~\ref{a10.110} concerning the definition of
rank) for generic $c_\la\in J$.  If $c_1$ is nondegenerate,
one can define the so-called recursion operator
$c_2\circ(c_1)^{-1}:TM\rightarrow TM$.  Its eigenvalues are
in involution with respect to $c_1$.  Off course, one
should impose additional conditions on $J$ in order that
these functions are independent and that they form a
"complete" set, i.e. the foliation defined by them is
lagrangian (not only coisotropic).

One can obtain examples of global symplectic
bihamiltonian structures considering holomorphic symplectic
manifolds $(M,\omega)$ and putting $c_1=\Re c,c_2=\Im c$,
where $c=(\omega)^{-1}$ is the holomorphic bivector field
inverse to  $\omega$. Off course, locally they are all the
same due to the Darboux theorem and they are quite not
interesting from the above point of view since the
recursion operator coincides with the complex structure and
has only constant eigenvalues $\pm\i$. However, these
bihamiltonian structures will be important for us and we
call them holomorphic symplectic.

The second class consists of degenerate bihamiltonian
structures, which are described by the condition
$max_\la\ra c_\la<\dim M$. Given such a structure, one can
construct the family of functions $F_0=
\linebreak \sum_{c\in J_0}Z_c$, where $J_0\subset J$ is a
subset of tensors of maximal rank and $Z_c$ stands for the
space of local Casimir functions of a Poisson tensor $c$.
It turns out that this family is in involution with respect
to any $c_\la \in J$. Again, keeping in mind the aim of
getting the completely integrable system, one should put
some restrictions on $J$.

An elegant and easily checkable condition which guarantees
that the family $F_0$ is locally "complete",
i.e.  defines a lagrangian foliation on a generic
symplectic leaf of $J$, is given by the Bolsinov-Brailov
theorem (see~\ref{a20.30}) and is formulated as follows:
\begin{description}
\item[($*$)] $\ra (\la_1c_1+\la_2c_2)(x)=R_0$ for any
$(\la_1,\la_2)\in\C^2\smallsetminus\{0\}$.
\end{description}
Here $R_0=\max_\la \ra c_\la$ and $x$ is a point in a
neighbourhood of which we are checking the "completeness"
of $F_0$.
Taking ($*$) as a starting point we define
complete bihamiltonian structures as those satisfying this
condition on an open dense subset.

The reader is referred to
papers~\cite{gz1}-\cite{gz5}  for the detailed exposition
of the geometric and algebraic aspects of bihamiltonian
structures based on the classical theory of pencils of
operators. Note that according to the terminology of these
articles the symplectic and complete bihamiltonian
structures mentioned above should be called (micro) Jordan
and (micro) Kronecker (cf.  Theorem~\ref{a20.65}, below).

The aim of this paper is to study some relations between
the above classes. More precisely, we study the reductions
of holomorphic symplectic bihamiltonian structures (by
means of real foliations) resulting in complete ones.
The construction inverse to such a reduction is called
a realization; hence the title of the paper.  Our main
theorem (see~\ref{main}) states the
completeness of the reduction $J'$ of the holomorphic
symplectic structure $J$ associated with the canonical
symplectic form $\omega$ on a generic coadjoint orbit
$M\subset\D$, where $\A$ is a complex Lie algebra from a
wide class including the semisimple algebras (cf.
Convention~\ref{d40.15}) and the reduction is performed
with respect to a real form $G_0\subset G$ of the (simple,
semisimple) complex Lie group adjoint to $\A$. Also, the
"first integrals", i.e.  the elements of the family $F_0$
corresponding to $J'$, are calculated
(Proposition~\ref{fin}). Let us make a few comments on the
proof.

Note that generic coadjoint $G_0$-orbits in $\D$ satisfy
the condition of $CR$-genericity (see
Definition~\ref{a10.150}).  The main theorem is a consequence
of a relatively simple criterion (Theorem~\ref{ns}) of
completeness for the reduction $J'$ of a holomorphic
symplectic structure $J$ by means of a real $CR$-generic
foliation ${\cal K}$ on $M$. We want to stress that the
assumption of the $CR$-genericity for ${\cal K}$ is natural
in the context discussed in Section~4. The study of
reductions without this assumption
seems reasonable, but more complicated.

In order to use the mentioned criterion to the
proof of the main theorem, one studies the auxiliary
complex Poisson pair $c,\tilde{c}$ on $\D$, where
$c:\D\rightarrow\D\wedge\D$ is the canonical linear Poisson
bivector and $\tilde{c}$ is its composition with the real
involution corresponding to the real form $\A_0\subset\A$.
It turns out that for the reduced Poisson pairs
$c',\bar{c}'$ and $c',\tilde{c}'$ the "first integrals"
coincide (see the proof of Theorem~\ref{main}). But it is
easy to calculate them for $c',\tilde{c}'$ using the method
similar to the classical method of argument translation
(\cite{mf}). This enable us to control rank of the
bivectors from $J'$ in the way required in the criterion.

The last two
essential ingredients of the proof are the
Gelfand-Zakharevich theorem about the structure of the pair
of bivectors in a vector space (see~\ref{a20.65}) and
some $CR$-geometric facts about the $G_0$-orbits in $\D$
(Section 6).

Note that the pair $c,\tilde{c}$ is defined canonically and
some of the results about it combined with that from Section
6 may be of independent interest (see Proposition~\ref{CR},
for example).

The paper is organized as follows. In Section~1 we recall
some definitions and facts from the theory of Poisson
manifolds and introduce the class of complex Poisson
structures, which are generalizations of the standard ones
to the case of the complexified tangent bundle. Holomorphic
Poisson structures are strictly contained in this class.
Also, we recall elementary  definitions from the theory of
$CR$-manifolds and adapt some of them to the symplectic
context.

Section~2 is devoted to bihamiltonian structures and their
relations with the completely integrable systems. We define
complete bihamiltonian structures, present some examples
and describe their structure from the point of view of the
Gelfand-Zakharevich theorem.

In Section~3 we prove that the
Poisson reduction $(c'_1,c'_2)$ of a Poisson pair $(c_1,c_2)$ is
again a Poisson pair under the requirement of the linear
independence for $c'_1,c'_2$. This result follows from the
natural behavior of the Schouten bracket with respect to
the reduction. We also study the relations between the
characteristic distributions of $c_\la$ and $c'_\la$.

The main theorem of Section~4 (see \ref{ns}) is the
criterion mentioned above. It is preceded by the
discussion of the linear algebraic aspects of the
reductions resulting in complete bihamiltonian
structures. In the end of this section a notion of minimal
realization is discussed.

The goal of Section~5 is to study the auxiliary Poisson pair
$c,\tilde{c}$ from the point of view of Section~2. In
particular, it is proved that it is complete and the
corresponding "first integrals" are calculated.

In section~6 we show that the coadjoint $G_0$-orbits are
$CR$-generic outside some $G_0$-invariant  real algebraic
set in $\D$ and calculate their dimension and
$CR$-dimension. We also show that they are isotropic with
respect to the canonical holomorphic symplectic form
$\omega$ on the corresponding $G$-orbit.

The concluding section is devoted to the formulation and
proof of the main theorem.

Off course, the inspiration for this paper is, besides the
mentioned papers of I.M.Gelfand and I.S.Zakharevich,
the theory of symplectic realizations for Poisson
structures (\cite{weinst}). Considering realizations of
degenerate bihamiltonian structures in symplectic
bihamiltonian structures different from holomorphic ones is
also meaningful (recently the author was informed by
Prof. F.J.Turiel that realizations in a symplectic bihamiltonian
structure of different kind, but also with constant
coefficients, give an elegant way for reconstructing the
bihamiltonian structure from its Veronese web,
cf.~\cite{gz2}). However, the author hopes that using of
the holomorphic structures opens a new perspective of
applying the complex-analytic methods to the theory of real
bihamiltonian structures.

We conclude this introduction by the following
conjecture: a generic real-analytic complete bihamiltonian
structure has a realization in a holomorphic symplectic
one; the double complex of differential operators related
to the problem of reconstruction the bihamiltonian
structure from its Veronese web (see~\cite{gz2}) is a kind
of reduction of the $d,d^c$-bicomplex.

\section{Complex Poisson structures and other
preliminaries}
\label{s10}

\abz\label{a10.10}\begin{abza}\rm
Let $M$ be a $\Ci$-manifold; write $\E (M)\ (\E^{\C} (M))$
for the space of $\Ci$-smooth real (complex) valued
functions on $M$.  We shall write $TM$ for the tangent
bundle and $T^\C M$ for its complexification.

All complex manifolds $M$ will be treated from the $\Ci$
point of view, so we shall not use special symbols for the
underlying real manifolds. The holomorphic tangent bundle
will be denoted by $T^{1,0}M$.

For a $\Ci$ vector bundle
$\pi:N\rightarrow M$, let $\G(N)$ denote the space of
$\Ci$-smooth sections of $\pi$.  Elements of $\G(\bigwedge
^2TM)\ (\G(\bigwedge ^2T^{\C} M))$ will be called (complex)
bivectors for short.
\end{abza}

\abz\label{a10.20}
\begin{defi}
A (complex) bivector $c\in \Gamma(\bigwedge^2TM)\ (\G(\bigwedge
^2T^{\C}M))$ is called
Poisson if $[c,c]=0$.
\end{defi}

Here
$[\,,]$ denotes the complex extension of the Schouten bracket
which associates a
trivector field $[c_1,c_2]\in \G(\bigwedge ^3T^{\C}M)$ to two
bivectors $c_1,c_2\in \Gamma(\bigwedge^2T^{\C}M)$.
The corresponding local coordinate
formula looks as follows:
\sabz\label{f10.20}
\begin{equation}
[c_1,c_2]^{ijk}(x)=\frac{1}{2}\sum_{\mbox{c.p.}
ijk}(c_1^{ir}(x)\frac{\d}{\d
x^r}c_2^{jk}(x)+c_2^{ir}(x)\frac{\d}{\d x^r}c_1^{jk}(x)),
\end{equation}
where $c_\alpha=c_\alpha^{ij}(x)\frac{\d}{\d
x^i}\wedge\frac{\d}{\d x^j}, \alpha=1,2,\
\sum_{\mbox{c.p.} ijk}$ denotes the  sum over the cyclic
permutations of $i,j,k$ and the summation convention over
repeated indices is used (the latter will be used
systematically in this paper).

\abz\label{10.15}
\begin{defi}
Let $M$ be a complex manifold. A holomorphic section of the
bundle $\bigwedge^2T^{1,0}M\subset T^\C M$ will be called
a holomorphic bivector. In particular, holomorphic
bivectors can be considered as complex ones and they will
be called holomorphic Poisson if, in addition, they are
Poisson in the sense of previous definition.
\end{defi}

\abz\label{a10.30}
\begin{defi}
A hamiltonian vector field $c(f)$ corresponding to a function
$f\in\E (M)\ (\E^{\C} (M))$ is obtained by the  contraction of the
differential $df$ and the Poisson bivector $c$ with respect to the first
index.
\end{defi}

\abz\label{a10.40}
\begin{prop}
A (complex) bivector $c$ is Poisson if and only if an operation
$\{\,,\}_c:\E (M)\times\E (M)\longrightarrow\E (M)\
(\E^{\C} (M)\times\E^{\C} (M)\longrightarrow\E^{\C} (M))$
given by
$$\{f,g\}_c=c(f)g,\ f,g\in\E (M)\ (\E^{\C} (M))$$
is a Lie algebra bracket over $\R\ (\C)$.

If $c$ is Poisson, then the map $f\mapsto c(f):\E
(M)\rightarrow \Vect(M)$, where $\Vect (M)$ is
a Lie algebra of smooth vector fields on $M$ with the
commutator bracket, is a Lie algebra homomorphism.
\end{prop}

\abz\label{a10.50}
\begin{defi}
$\{\,,\}_c$ is called the Poisson bracket corresponding to
a (complex) Poisson bivector $c$. A family of functions
$F\subset\E (M)\ (\E^{\C} (M))$ is involutive with
respect to $c$ if $\{f,g\}_c=0$ for each two functions
$f,g\in F$.
\end{defi}

\abz\label{a10.100}
\begin{defi}
Consider a (complex) bivector $c$ at $x\in M$ as a map
$c^\sharp _x: T_x^*M\longrightarrow T_xM\
((T_x^{\C}M)^*\longrightarrow T_x^{\C}M)$ evaluating the
first argument of the bivector on a $1$-form.  Kernel $\ke
c(x)$ and rank $\ra c(x)$ of $c$ at $x$ are defined as that
of $c^\sharp_x$. We say that $c$ is nondegenerate if
$c^\sharp$ is an isomorphism.
A complex bivector of type $(2,0)$ on a complex manifold
$M$ is called nondegenerate in the holomorphic sense if the
restricted sharp map $c^\sharp :(T^{1,0}M)^*\longrightarrow
T^{1,0}M$ is an isomorphism.
\end{defi}

\abz\label{a10.60}
\begin{defi}
A characteristic subspace $P_{c,x}$ of a (complex) bivector
$c$ at $x\in M$ is defined as $\im c_x^\sharp$.
A generalized distribution of subspaces
$P_c\subset TM\ (T^{\C}M)$ is said to be a characteristic
distribution for the  bivector $c$.  \end{defi}

Note that a complex
bivector of type $(2,0)$ nondegenerate in holomorphic sense
is not nondegenerate since $P_{c,x}=T^{1,0}M\not =
T^\C M$.  We shall usually understand the nondegeneracy of
holomorphic bivectors in the holomorphic sense.

\abz\label{a10.70}
\begin{th}
(\cite{kir}) Let $c$ be a real Poisson bivector.
The generalized distribution $P_c$ is
completely integrable, i.e. there exists a tangent to $P_c$
generalized foliation $\{S_\alpha\}_{\alpha\in I}$ on $M$:
$T_xS_\alpha=P_{c,x}$ for any $\alpha\in I$ and for any $x\in
S_\alpha$. The restriction of $c$ to each $S_\alpha$ is a
nondegenerate Poisson bivector; consequently, $S_\alpha$
are symplectic manifolds with the symplectic forms
$\omega_\alpha=(c|_{S_\alpha})^{-1}$.  \end{th}

Here and subsequently the 2-form $\omega$ inverse to a
nondegenerate bivector $c$
is defined as follows. If $\wedge^2 c^\sharp$ is the
extension of the sharp map defined above to the second
exterior power of $T^*M$, then $\omega=(\wedge^2
c^\sharp)^{-1}(c)$. The inverse to a nondegenerate 2-form
bivector is defined similarly.

The above theorem is also true in the complex analytic
category if we understand $P_c$ as a holomorphic subbundle
in $T^{1,0}M$ and the nondegeneracy in the holomorphic
sense.  The definition of inverse objects in this case is
analogous to real one.

\abz\label{a10.80}
\begin{defi}
The submanifolds
$S_\alpha$ are called symplectic leaves of a Poisson
bivector $c$.
\end{defi}

\abz\label{a10.90}
\begin{prop}
Given  a complex Poisson bivector $c\in\Gamma(\bigwedge^2T^{\C}M)$,
its characteristic distribution $P_c\subset T^{\C}M$ is
involutive, i.e.
$$[v,w](x)\in P_{c,x}$$
for any complex valued vector fields $v,w$ such that
$v(x),w(x)\in P_{c,x}, x\in M$.
\end{prop}

In general, one can say nothing about the complete integrability
of $P_c$ even if one understands this in spirit of the
Newlander-Nierenberg theorem. A nonconstant rank of the
subspaces $P_{c,x}$ or $P_{c,x}\bigcap\overline P_{c,x}$ (the
overline means the complex conjugation)  may be the obstruction here as
well as some other reasons (see \cite{treves}).

\abz\label{dex}
\begin{exa}\rm
Let $M=\C^3$ with coordinates $z_1,z_2,z_3$,
$c=\bar{z}_1\frac{\partial }{\partial
z_2}\wedge\frac{\partial }{\partial
z_3}+\bar{z}_2\frac{\partial }{\partial
z_3}\wedge\frac{\partial }{\partial
z_1}+\bar{z}_3\frac{\partial }{\partial
z_1}\wedge\frac{\partial }{\partial z_2}$,
$f=|z_1|^2+|z_2|^2+|z_3|^2$. The bivector $c$ is obviously
Poisson. Since $c(f)=0$, its characteristic subspace
$P_{c,x}$ is equal to the $(1,0)$-tangent space
(cf.~\ref{a10.145}) to the $5$-dimensional sphere
$S\subset M$ centered in $0$ and passing through $x$.
Off course, this example is related to the Lie algebra
$so(3)$. We shall generalize it in Section 5.
\end{exa}

\abz\label{a10.110}
\begin{conv}\rm
In the sequel, all Poisson bivectors will
be assumed to have
maximal rank on an open
dense subset in $M$. For real Poisson bivectors this is
equivalent to the following: the union of symplectic leaves
of maximal dimension is dense in $M$.
\end{conv}
\abz\label{a10.120}
\begin{defi}
Let $c$ be a (complex) Poisson bivector. Define $\ra c$ as
$\max_{x\in M}\ra c(x)$.
\end{defi}

\abz\label{a10.130}
\begin{defi}
A Casimir function $f\in\E (U)\ (\E^{\C} (U))$ over an open set $U\subset
M$ for a (complex) Poisson bivector $c$ is defined by the condition $c(f)=0$.
A space of all Casimir functions  over $U$ for $c$ is denoted by
$Z_c(U)$ and $Z_{c,x}$ stands for the space of germs of Casimir functions
at $x,x\in M$.
\end{defi}

Note that if $c$ is real and $\ra c<\dim M$ there exist local nontrivial
Casimir functions and their differentials at $x$ span
$\ke c(x)$, provided that $x$ is taken from a symplectic leaf of
maximal dimension. This is not true concerning the global
Casimir functions: it is easy to construct a Poisson bivector
$c$ with $\ra c<\dim M$ possessing only trivial ones.

\abz\label{a10.140}
\begin{defi}
Let $(M,\omega),\,\dim M=2n$, be a symplectic manifold. A
submanifold $L\subset M$ is called
\begin{itemize}
\item coisotropic if $(T_xL)^{\bot\omega(x)}\subset T_xL$ for
any $x\in L$;
\item isotropic if $(T_xL)^{\bot\omega(x)}\supset T_xL$ for any
$x\in L$;
\item lagrangian if $(T_xL)^{\bot\omega(x)}=T_ xL$ for any
$x\in L$.
\end{itemize}
A foliation ${\cal L}$ on $M$ is coisotropic (isotropic,
lagrangian) if so is its every leaf.
\end{defi}

Here $\bot\omega(x)$ stands for a skew-orthogonal complement in
$T_xM$ with respect to $\omega(x)$. For the third case the
following definition is equivalent: $\dim L=n$ and
$\omega|_{TL}\equiv 0$.

\noindent\abz\label{a10.145}\begin{abza}\rm
We shall need a specific generalization of this definition in
the complex case. Let $M$ be a complex manifold with the
complex structure ${\cal J}:TM\longrightarrow TM$. Consider
a $\Ci$-smooth submanifold $L\subset M$.  Write $T_x^{CR}L$
for $T_xL\bigcap{\cal J}T_xL$ and $T_x^{1,0}L$ for
$T_x^{\C}L\bigcap T_x^{1,0}M,\, x\in L$. Another definition for
$T_x^{1,0}L$ is the following: $T_x^{1,0}L=\{v-\i {\cal J}v;v\in
T_x^{CR}L\}$.
\end{abza}

\abz\label{a10.150}
\begin{defi}
(\cite{treves},\cite{boggess}) $L$ is called a $CR$-submanifold in
$M$ if \linebreak $\dim T_x^{1,0}L$ is constant along $L$;
we say that this number is
$CR$-dimension of $L$; L is  generic  (completely real) if
$T_xL+{\cal J}T_xL=T_xM$ ($T_xL\oplus{\cal J}T_xL=T_xM$) for each $x\in L$.
\end{defi}

If a generic $CR$-submanifold $L$ is given by the equations
$\{f_1=\alpha_1,\ldots,f_k=\alpha_k\},f_i\in\E (M)$, such that
$df_1\wedge\ldots\wedge
df_k\neq 0$ along $L$, $\alpha_i\in \R$, then
$\d f_1\wedge\ldots\wedge\d f_k\neq 0$ along $L$ and
$T_x^{1,0}L=\{\d f_1(x),\ldots ,\d f_k(x)\}^{\perp 1,0}$, where $\d$ is
the (1,0)-differential on $M$, $\perp 1,0$ denotes the
annihilator in $T_x^{1,0}M$.

\abz\label{a10.155}
\begin{defi}
A foliation
${\cal L}$ on $M$ is a generic (completely real) $CR$-foliation if its each
leaf
is a generic (completely real) $CR$-submani-\lbr fold.
\end{defi}

\abz\label{a10.160}
\begin{defi}
Let $(M,\omega)$ be a holomorphic symplectic
manifold. A $CR$-submanifold $L\subset M$ is
\begin{itemize}
\item $CR$-coisotropic if
$(T_x^{1,0}L)^{\bot\omega(x)}\subset T_x^{1,0}L$ for any
$x\in L$; \item $CR$-isotropic if
$(T_x^{1,0}L)^{\bot\omega(x)}\supset T_x^{1,0}L$ for any
$x\in L$; \item $CR$-lagrangian if
$(T_x^{1,0}L)^{\bot\omega(x)}= T_x^{1,0}L$ for any $x\in
L$.
\end{itemize}
A $CR$-foliation ${\cal L}$ on $M$ is
said to be $CR$-coisotropic ($CR$-isotropic,
$CR$-lagrangian) if so is its every leaf.
\end{defi}

Here $\bot\omega(x)$ denotes a skew-orthogonal complement in
$T_x^{1,0}M$ with respect to the (2,0)-form $\omega(x)$.

Suppose ${\cal L}$ is generic and consists of
the common level sets of the functions $f_1,\ldots,f_k\in\E (M),
k\leq n=(1/2)\dim_\C M$. Then ${\cal L}$ is
$CR$-coisotropic if and only if the family
$\{f_1,\ldots,f_k\}$ is involutive with respect to the
holomorphic Poisson bivector $c=(\omega)^{-1}$. In
particular, if  $k=n$ one gets $CR$-lagrangian
foliation.

\abz\label{a10.170}
\begin{defi}
Let $(M,\omega), \dim M=2n$, be a symplectic manifold.
A completely integrable system  on $M$ is
defined as a family of functions ${\cal F}\subset\E (M)$
involutive with respect to $c=(\omega)^{-1}$ and containing a subfamily of
$n$ functions that are functionally independent almost
everywhere on $M$. In other words, a completely integrable
system on $M$ is a lagrangian foliation ${\cal L}$ on an open
dense subset in $M$.
\end{defi}

We conclude this section by recalling main definitions
concerning hamiltonian actions of Lie groups
(see~\cite{weinst2} for details).

\abz\label{10.17}
\begin{defi}
Let $G$ be a connected Lie group with the Lie algebra $\A$.
Assume it is acting on a Poisson manifold $M$ with the
Poisson bivector $c$, i.e. a Lie algebra
homomorphism $\rho:\A\rightarrow \Vect (M)$ is given.   The
action is called hamiltonian if there exists a Lie algebra
homomorphism $\psi:\A\rightarrow \E (M)$ such that the
following diagram is commutative
$$
\begin{array}{ccc}
\A  & \stackrel{\psi}{\longrightarrow}& \E(M) \\
\parallel &  & \downarrow  c(\cdot) \\
\A &\stackrel{\rho}{\longrightarrow}
&\Vect (M),
\end{array}
$$
where $c(\cdot)$ is a Lie algebra homomorphism of taking
the hamiltonian vector field (see
Proposition~\ref{a10.40}).

%The map $x\mapsto \varphi_x:M\rightarrow \D$ defined by
%$\varphi_x(v)=\psi(v)(x), v\in\A$ is called the moment
%map.
\end{defi}

\section{Bihamiltonian structures and completeness}
\label{s2}

Let $M$ be a $\Ci$-manifold.

\abz\label{a10.180}
\begin{defi}
Two linearly independent (complex) Poisson bivectors $c_1,c_2$
on $M$ form a (complex) Poisson pair if
$c_\la=\la_1c_1+\la_2c_2$ is  a Poisson bivector
for any $\la=(\la_1,\la_2)\in\R^2$ ($\C^2$).  \end{defi}

\abz\label{a10.190}
\begin{prop}
A pair of linearly independent (complex) bivectors $(c_1,c_2)$ is
Poisson if and only if $[c_1,c_1]=0, [c_1,c_2]=0,
[c_2,c_2]=0$.
\end{prop}

%\abz\label{a10.200}
%\begin{prop}
%Let $(c_1,c_2)$ be a pair of linearly independent real bivectors.
%Then $(c_1,c_2)$ is a Poisson pair if and only if  so is the
%pair $(c,\bar{c})$, where $c=c_1+\i c_2, \bar{c}=c_1-\i c_2$.
%\end{prop}

\abz\label{a10.210}
\begin{defi}
Let $M$ be a $\Ci$-manifold. A (complex) bihamiltonian structure
on $M$ is defined as a two-dimensional linear subspace
$J=\{c_\la\}_{\la\in{\cal S}}$ of (complex) Poisson bivectors on $M$
parametrized by a two-dimensional vector space ${\cal S}$ over
$\R$ ($\C$). The trivial bihamiltonian structure is the
zero-dimensional linear subspace in $\G(\bigwedge^2TM)$.
\end{defi}

It is clear that every Poisson pair generates a
nontrivial bihamiltonian structure and the transition from
the latter one  to a Poisson pair corresponds to a choice
of basis in ${\cal S}$. We shall write $(J,c_1,c_2)$ for a
bihamiltonian structure $J$ with a chosen Poisson pair
$(c_1,c_2)$ generating $J$.

\abz\label{a10.220}
\begin{defi}
Let $(J,c_1,c_2)$ be a bihamiltonian structure.
A complex bihamiltonian structure
$$J^\C=\{\la_1c_1+\la_2c_2;(\la_1,\la_2)\in\C^2\}$$
is called the complexification of $J$.
\end{defi}

\abz\label{a10.225}
\begin{prop}
A complex bihamiltonian structure $J$ is the  complexification of
some real one if and only if one can choose a generating $J$
Poisson pair $(c,\bar{c})$, where $c\in \bigwedge^2(T^{\C}M)$,
the bar stands for the complex conjugation.
\end{prop}

\noindent{\sc Proof.} Generate $J'$ by $\Re c, \Im c$. Then
$J=(J')^\C$.  Conversely, one checks that for
$(J',c_1,c_2)$ the complexification $(J')^\C$ is generated
by $c_1\pm\i c_2$. q.e.d.

\abz\label{a10.230}
\begin{defi}
Let $J$ be a (complex) bihamiltonian structure and let
$J_0\subset J$ be a subfamily of (complex) Poisson bivectors of
maximal rank $R_0$ (the set $J\smallsetminus J_0$ is at most a finite sum of
1-dimensional subspaces). We say that
$J$ is symplectic if $\ra c_\la=\dim M$ for any $c_\la\in
J_0$ and that $J$ is degenerate otherwise.
\end{defi}

\noindent\abz\label{a10.240}
\begin{exa}\rm
Consider a family $J^\C$
generated by a pair $(c,\bar{c})$, where $c=(\omega)^{-1}$
is a complex Poisson bivector inverse to a holomorphic
symplectic form $\omega$ on a complex symplectic manifold
$M$. Since $c$ is holomorphic and $\bar{c}$ is
antiholomorphic, we have $[c,\bar{c}]=0$.  Thus $J^\C$ is a
bihamiltonian structure. By Proposition \ref{a10.225} it is
the complexification of the real bihamiltonian structure
$(J,\Re c, \Im c)$. This example
is fundamental for the
paper and we shall need the following fact.
\end{exa}

\abz\label{a10.245}
\begin{prop} Let $M, \omega, c$ and $J^\C$ be as in Example
\ref{a10.240}.  Then $J^\C$ is  symplectic and the only
degenerate bivectors in $J^\C$ are those proportional to
$c$ and $\bar{c}$. Moreover, $\ra_\C
c=\ra_\C\bar{c}=\frac{1}{2}\dim_{\C}M, P_c=T^{1,0}M,
P_{\bar{c}}=T^{0,1}M$.  \end{prop}

\noindent{\sc Proof.} The last assertion is obvious as well as the
following equality
$$c\con \overline{\omega}=0,$$
where $\con $ stands for the contraction with respect to the first index.
For $\omega_1=\Re\omega,\omega_2=\Im\omega,c_1=\Re c,c_2=\Im c$ this
implies
$$c_1\con \omega_1+c_2\con \omega_2=0,c_2\con \omega_1-c_1\con \omega_2=0.$$
We have for $\la\in\C$
\begin{eqnarray*}
(c_1+\la c_2)\con (\omega_1-\la\omega_2)=
c_1\con \omega_1-\la^2c_2\con \omega_2-
\la(c_1\con \omega_2-c_2\con \omega_1)=\\(1+\la^2)c_1\con \omega_1=(1+\la^2)
\frac{1}{4}id_{TM}.
\end{eqnarray*}
The last equality is verified directly in the Darboux local
coordinates. Thus $(c_1+\la c_2)^{-1}=\frac{4}{1+\la^2}(\omega_1-\la\omega_2)$.
Since $c_1\con \omega_1=-c_2\con\omega_2$, $c_2$ is also nondegenerate. q.e.d.

\abz\label{a10.247}
\begin{defi}
Let $(M,\omega)$ be a complex symplectic manifold. The
bihamiltonian structure $J$ and its complexification $J^\C$
from Example \ref{a10.240} are  called  holomorphic
symplectic.  \end{defi}

\abz\label{a10.248}\begin{abza}\rm
Given a (complex) bihamiltonian structure $J$, let $F_0$ denote the
space $\sum_{c\in J_0}Z_c(M))\subset\E (M)$.
\end{abza}

The following theorem shows how  the degenerate
bihamiltonian structures can be applied for constructing the completely
integrable systems.

\abz\label{a10.250}
\begin{th}
Let $J$ be a degenerate (complex) bihamiltonian structure on $M$. A family
$F_0$ is involutive
with respect to any $c_\la\in J$.
\end{th}

\noindent{\sc Proof.} Let $c_1,c_2\in J_0$ be linearly independent,
$f_i\in Z_{c_i}, i=1,2$. Then
\sabz\label{f1}
\begin{equation}
\{f_1,f_2\}_{c_\la}=(\la_1c_1(f_1)+\la_2c_2(f_1))f_2=-\la_2c_2(f_2)f_1=0.
\end{equation}
Now it remains to prove that for any $c\in J_0,f_i\in Z_c,
i=1,2$, one has
$\{f_1,f_2\}_{c_\la}=0$.
For that purpose we first rewrite
(\ref{f1}) as
\sabz\label{f2}
\begin{equation}
c_\la(x)(\phi_1,\phi_2)=0,
\end{equation}
where $\phi_i\in \ke c_i(x),\ i=1,2,\ x\in M$, and the
lefthandside denotes a contraction of the bivector with two covectors.
Second, we fix $x$ such that $\ra c(x)=R_0$
and approximate $df_2|_x$ by a sequence of elements
$\{\phi^i\}_{i=1}^\infty,\ \phi^i\in \ke c^i(x)$, where $c^i\in
J_0,\ i=1,2,\dots,$ is linearly independent with $c$. Finally,
by (\ref{f2}) we get $c_\la(x)(df_1|_x,\phi^i)=0$ and by the
continuity $\{f_1,f_2\}_{c_\la}(x)=0$. Since the set of such
points $x$ is dense in $M$, the proof is finished.
q.e.d.

In fact this theorem is true for the local Casimir
functions (for the germs of Casimir functions).

\abz\label{20.05}
\begin{defi}
The functions from the family $F_0$ (see
\ref{a10.248}) are called  (global) first integrals of
the bihamiltonian structure $J$.  The family of functions
$\sum_{c\in J_0}Z_c(U)$
($\sum_{c\in J_0}Z_{c,x}$) is denoted by
$F_0(U)$ ($F_{0,x}$) and its elements are
called local first integrals over an open
$U\subset M$ (germs of first integrals at $x\in
M$).  \end{defi}

In order to obtain a completely integrable
system from Casimir functions one should
require additional assumptions on the bihamiltonian structure $J$.
Off course, the condition of completeness
given below concerns the local Casimir functions (in fact their germs)
and may be insufficient for
obtaining the completely integrable system. However, it is of
use if the local Casimir functions are restrictions of the
global ones (see Example \ref{a20.50}, below).

Given a characteristic distribution $P_c\subset TM$ ($T^{\C}M$) of some
(complex) Poisson bivector and a point $x\in M$, let $P_{c,x}^*$ denote a
dual space to $P_{c,x}$. Any functional $\phi\in T_x^*M$ ($(T_x^{\C}M)^*$)
can be regarded as an element of $P_{c,x}^*$ called the restriction of $\phi$
to $P_{c,x}$.

\abz\label{a20.10}
\begin{defi}
(\cite{bols}) Let $J$ be a (complex) bihamiltonian structure;
fix some $c_\la\in J$.

$J$ is called complete at a point $x\in M$ with
respect to $c_\la$ if a
linear subspace of $P_{c_\la,x}^*$ generated over  $\R$ ($\C$) by the
differentials of
the germs $f\in F_{0,x}$ restricted to $P_{c_\la,x}$ has dimension
$\frac{1}{2}\dim_{\R} P_{c_\la,x}$ ($\frac{1}{2}\dim_{\C}P_{c_\la,x}$).
\end{defi}

\abz\label{a20.20}
\begin{prop}
A (complex) bihamiltonian structure
$J$ is complete with respect to $c_\la\in J$ at a point $x\in M$
if and only
if $\dim(\bigcap_{c\in J_0}\! P_{c,x})\linebreak=\frac{1}{2}\dim P_{c_\la,x}$.
\end{prop}

\noindent {\sc Proof} is obvious.

The following theorem is due to A.Brailov (see \cite{bols}, Theorem 1.1 and
Remark after it).

\abz\label{a20.30}
\begin{th}
A (complex) bihamiltonian structure
$J$ is complete with respect to $c_\la\in J_0$ at a point $x\in M$
such that $P_{c_\la,x}$ is of maximal dimension if and only
if the following condition holds
\begin{description}
\item[$(*)$] $\ra c(x)=R_0$ for any
$c\in J^\C\smallsetminus \{0\}$ ($J\smallsetminus\{0\}$),
\end{description}
where $R_0$ is as in \ref{a10.230}.
\end{th}

\noindent{\sc Proof} of this theorem is a consequence of
the following linear algebraic fact.

\abz\label{20.35}
\begin{prop}
(\cite{bols}) Let $V$ be a vector space over $\R$ ($\C$)
and let $J$  be a two
dimensional linear subspace in $\bigwedge^2V$.
In the real case we let $J^\C\subset
\bigwedge^2V^\C$ denote the complexification of the
subspace $J$.  We write $J_0\subset J$ for the subset of
bivectors of maximal rank $R_0$ and $F_0\subset V^*$ for
the subspace generated by the kernels of bivectors from
$J_0$.  Let $c^\sharp :V^*\rightarrow V$ stand for the
corresponding sharp map of $c\in \bigwedge^2V$
(cf.~\ref{a10.100}).

Then, given a bivector $c_\la\in J_0$, the following two
conditions are equivalent:
\begin{description}
\item[($i$)] $\dim (F_0|_{P_\la})
=1/2\dim P_\la$, where $F_0|_{P_\la}=
\Span\{\{f|_{P_\la}\}_{f\in F_0}\}\subset P_\la^*$ and
$P_{\la}=c_\la^\sharp (V^*)$;
\item[($ii$)] $\ra c=R_0$ for
any $c\in J^\C\smallsetminus\{ 0\}$ ($J\smallsetminus\{ 0\}$).
\end{description}
\end{prop}

\noindent{\sc Proof.}
We reproduce the proof from~\cite{bols} with a small
completion.

We perform the proof in the following four steps.

First, we observe that for any two bivectors
$a,b\in J\smallsetminus\{0\}$ one has the equality $a^\sharp
(F_0) =b^\sharp (F_0)$. Indeed, suppose that $a,b$ are
linearly independent. The subspace $F_0$ is generated by a
finite number of kernels $\ker b_1,\ldots ,\ker b_s, b_i\in
J_0$. Without loss of generality, we may assume that
$b_i=\alpha_ia+\beta_ib$, where $\alpha_i,\beta_i\not =
0$. Since $(\alpha_ia^\sharp+\beta_ib^\sharp)(\ker
b_i)=0, i=1,\ldots ,s$, then $a^\sharp(\ker
b_i)=b^\sharp(\ker b_i)$ and, consequently,
$a^\sharp (F_0)=b^\sharp(F_0)$.

Second,  consider the skew-orthogonal
complement $\tilde{F}_0=(F_0)^{\bot b}=(b^\sharp
(F_0))^\bot$ and note that: 1) it does not depend on $b\in
J\smallsetminus\{ 0\}$ (previous step); 2)
$F_0\subset\tilde{F}_0$ (the skew-orthogonal complement of
any subspace in $V^*$ with respect to any $b\in
J\smallsetminus\{0\}$ contains $\ker b$, in particular
$\tilde{F}_0\supset\ker b, b\in J\smallsetminus \{0\}$); 3) if
$a\in J_0$, then $b^\sharp (\tilde{F}_0)\subset a^\sharp
(\tilde{F}_0)$ for any $b\in J$ (this is equivalent to
$(F_0^{\bot b})^{\bot b}\supset (F_0^{\bot a})^{\bot a}$ or
$F_0+\ker b\supset F_0+\ker a=F_0$).

Third, given two linearly independent bivectors $a,b\in
J$, with \linebreak $\ra a=R_0$, we define a "recursion"
operator $\Phi:\tilde{F}_0/F_0\rightarrow \tilde{F}_0/F_0$
by the formula $\Phi (\pi
(\xi))=\pi((a^\sharp)^{-1}b^\sharp(\xi))$, where $\xi\in
\tilde{F}_0$ and $\pi:\tilde{F}_0\rightarrow
\tilde{F}_0/F_0$ is the natural projection. The operator is
correctly defined due to the conditions $a(F_0)=b(F_0)$,
$b(\tilde{F}_0)\subset a(\tilde{F}_0)$, and $\ker a\subset
F_0$. It is easy to see that the eigenvalues of $\Phi$ are
precisely those $\la\in\C$ for which $\ra (a-\la b)<R_0$.
In particular ($ii$) holds if and only if $\Phi$ does not
have eigenvalues, i.e. $F_0=\tilde{F}_0$.

Finally, we use the following sequence of subspaces and
relations between them
\begin{eqnarray*}
F_0=\pi^{-1}(\pi (F_0))=\pi^{-1}(F_0|_{P_\la})\subset
\pi^{-1}((c^\sharp (F_0|_{P_\la}))^{\bot
_\la})= \\ (c_\la^\sharp(F_0|_{P_\la}))^\bot=
(c_\la^\sharp(F_0))^\bot= \tilde{F}_0,
\end{eqnarray*}
where $\pi::V^*\rightarrow V/P_\la^\bot\cong P_\la^*$ is the
canonical projection and $\bot_\la$ is the annihilator in
the sense of the dual pair $(P_\la,P_\la^*)$. The essential
moment here is that $\ker\pi=\ker c_\la\subset F_0$; this
implies the first equality. The only inclusion in this
sequence is the equality, i.e.  $F_0|_{P_\la}$ is a
lagrangian subspace with respect to $C_\la|_{P_\la}$, if
and only if condition ($i$) holds. q.e.d.

Theorem~\ref{a20.30} shows that $J$ is complete with
respect to a fixed $c_\la\in J_0$ at a point $x$ such that
the dimension $P_{c_\la,x}$ is maximal if and only if
$J=J_0\bigcup\{0\}$ and $J$ is complete at $x$ with respect
to any nontrivial $c_\la\in J$. This motivates the next
definition.

\abz\label{a20.40}
\begin{defi}
Let $(J,c_1,c_2)$ be a (complex) bihamiltonian structure.
The structure $J$ (the pair $(c_1,c_2)$)
is complete at a point $x\in M$ if condition
$(*)$ of Theorem \ref{a20.30} holds at $x$. $J$ ($(c_1,c_2)$)
is called complete if it is so at any point from some
open and dense subset in $M$. The trivial bihamiltonian
structure is complete by definition.
\end{defi}

\abz\label{d20.45}
\begin{cor}
Let a bihamiltonian structure $J$ be complete at any $x$
from some sufficiently small open set $U$. Then the
functions from $F_0(U)$ (see~\ref{a10.248}) define a
foliation ${\cal L}$ on $U$ that is lagrangian in any
symplectic leaf $S_\la$ of any $c_\la|_U, c_\la\in
J|setminus\{0\}$. On the overlap of two such sets the
corresponding foliations coincide.
\end{cor}

\noindent{\sc Proof.} The first assertion follows from
Theorem~\ref{a20.30}. The second one is a consequence of
the uniqueness of the set of local first integrals of a
degenerate bihamiltonian structure. q.e.d.

\abz\label{d20.47}
\begin{defi}
The foliation ${\cal L}$ described in
Proposition~\ref{d20.45} will be called the bilagrangian
foliation of a complete bihamiltonian structure.
\end{defi}

\noindent\abz\label{a20.50}\begin{exa}\rm

(Method of argument translation, see \cite{mf},
\cite{bols}.) Let $\A$ be a nonabelian Lie algebra, $\D$
its dual space. Fix a basis $\{e_1,\ldots,e_n\}$ in $\A$
with the structure constants $\{c_{ij}^k\}$.
The standard linear Poisson bivector on $\D$ is defined as
$$c_1(x)=c_{ij}^kx_k\frac{\d}{\d x_i}\wedge\frac{\d}{\d x_j},$$
where $\{x_k\}$ are linear coordinates in $\D$
corresponding to $\{e_1,\ldots,e_n\}$. In more invariant
terms $c_1$ is described as dual to the Lie-multiplication
map $[\,,]:  \A\wedge\A\longrightarrow \A$. It is
well-known that the symplectic leaves of $c_1$ are the
coadjoint orbits in $\D$. Now define $c_2$ as a bivector
with constant coefficients $c_2=c(a)$, where $a$ is a fixed
point on any leaf of maximal dimension. It turns out that
$c_1,c_2$ form a Poisson pair and it is easy to describe
the set $I$ of points $x$ for which condition $(*)$ fails.
Consider the complexification $(\D)^\C\cong (\A^\C)^*$ and
the sum $\Sing(\A^\C )^*$ of symplectic leaves of
nonmaximal dimension for the complex linear bivector
$c_{ij}^kz_k\frac{\d}{\d z_i}\wedge\frac{\d}{\d z_j},$
where $z_j=x_j+\i y_j,\ j=1,\ldots,n,$ are the corresponding
complex coordinates in $(\D)^\C$. Then $I$ is equal
to the intersection of the sets $\A^*\subset(\D)^\C$ and
$\overline{a,\Sing(\A^\C )^*}$,
where $\overline{a,\Sing(\A^\C )^*}$ denotes a cone of
complex $2$-dimensional subspaces passing through $a$ and
$\Sing(\A^\C )^*$.
In particular, $(c_1,c_2)$ is complete for a semisimple
$\A$ since $\Sing\,\A^{\C*}$ has codimension at least $3$.
Note that this gives rise to completely integrable systems,
since the local Casimir functions on $\D$ are restrictions
of the global ones, i.e.  the invariants of the coadjoint
action.  \end{exa}

\noindent\abz\label{a20.60}\begin{exa}\rm
(Bihamiltonian structure of general position on
an odd-dimensional manifold, see \cite{gz2}.) Consider a pair of
bivectors $(a_1,a_2)$, $a_i\in\bigwedge^2V,  i=1,2$, where $V$ is a
$(2m+1)$-dimensional vector space; $(a_1,a_2)$ is in general
position if and only if is represented by the Kronecker block of
dimension $2m+1$, i.e.
\sabz\label{krone}
\begin{equation}
\begin{array}{l}
a_1=p_1\wedge q_1+p_2\wedge q_2+\cdots+p_m\wedge q_m\\
a_2=p_1\wedge q_2+p_2\wedge q_3+\cdots+p_m\wedge q_{m+1}
\end{array}
\end{equation}
in an appropriate basis $p_1,\ldots p_m,q_1,\ldots,
q_{m+1}$ of $V$.  A bihamiltonian structure $J$ on a
$(2m+1)$-dimensional $M$ is in general position if and only if
the pair $(c_1(x),c_2(x))$ is so for any $x\in M$. Such $J$ is
complete: it is easy to prove that $J=J_0\bigcup\{0\},\ \dim\bigcap_{c\in
J}P_c(x)=n$ and then use Proposition \ref{a20.20}. In general, a complete
Poisson pair at a point is the direct sum of the Kronecker
blocks and the zero pair as the corollary of the next theorem shows.
This theorem is a reformulation of the classification result for
pairs of $2$-forms in a vector space
(\cite{gz1}, \cite{gz3}).
\end{exa}

\abz\label{a20.65}
\begin{th}
Given a finite-dimensional vector
space $V$ over $\C$ and a pair of bivectors $(c_1,c_2),\
c_i\in\bigwedge^2V,$ there exists a direct decomposition
$V=\oplus_{j=1}^kV_j,\ c_i=\sum_{j=1}^kc_i^{(j)},
\ c_i^{(j)}\in\bigwedge^2V_j,\ i=1,2,$ such that each pair
$(c_1^{(j)},c_2^{(j)})$ is from the following list:
\begin{description}
\item[(a)] the Jordan block: $\dim V_j=2n_j$ and in an
appropriate basis of $V_j$ the matrix of $c_i^{(j)}$ is
equal to $$\left(\begin{array}{cc} 0&A_i\\ -A_i^T&0
\end{array}\right),i=1,2,$$
where $A_1=I_{n_j}$ (the unity $n_j\times n_j$-matrix) and
$A_2=J_{n_j}^\la$ (the Jordan block with the eigenvalue $\la$);
\item[(b)] the Kronecker block: $\dim V_j=2n_j+1$ and in an
appropriate basis of $V_j$ the matrix of $c_i^{(j)}$ is
equal to $$\left(\begin{array}{cc} 0&B_i\\ -B_i^T&0
\end{array}\right), i=1,2,$$
where
${\renewcommand{\arraycolsep}{-.005cm}
B_1=\left(\begin{array}{cccccc}
1&0&0&\ldots&0&0\\
0&1&0&\ldots&0&0\\
 & & &\ldots& & \\
0&0&0&\ldots&1&0
\end{array} \right),\ B_2=\left(\begin{array}{cccccc}
0&1&0&\ldots&0&0\\
0&0&1&\ldots&0&0\\
 & & &\ldots& & \\
0&0&0&\ldots&0&1
\end{array}\right)}$ ($(n_j+1)\times n_j$-matrices).
\item[(c)] the trivial Kronecker block: $\dim V_j=1$,
$c_1^{(j)}=c_2^{(j)}=0$;
\end{description}
\end{th}

\abz\label{a20.67}
\begin{cor}
Let $J$ be a (complex) bihamiltonian structure. It is complete at a point
$x\in M$ if and only if a pair $(c_1(x),c_2(x)),\
c_i(x)\in\bigwedge^2(T_x^{\C}M), i=1,2,$ does not contain the
Jordan blocks in its decomposition.
\end{cor}

\noindent {\sc Proof} follows from the definition of completeness.

The following example of a complete Poisson pair shows that
the structure of decomposition to the Kronecker blocks may
change from point to point.

\abz\label{20.68}
\begin{exa}\rm (\cite{panas}) Let $M=\R^6$ with coordinates
$(p_1,p_2,q_1,\ldots,q_4)$, $c_1=\frac{\d}{\d
p_1}\wedge\frac{\d}{\d q_1}+\frac{\d}{\d
p_2}\wedge\frac{\d}{\d q_2}, c_2=\frac{\d}{\d
p_1}\wedge(\frac{\d}{\d q_2}+q_1\frac{\d}{\d q_3})+\frac{\d}{\d
p_2}\wedge\frac{\d}{\d q_4}$. Here we have: two $3$-dimensional
Kronecker blocks on $M\smallsetminus H, H=\{q_1=0\}$; the
$5$-dimensional Kronecker block and the $1$-dimensional zero
block on the hyperplane $H$.
\end{exa}

\abz\label{20.69}
\begin{rem}\rm The decomposition of a pair of bivectors
$c_1,c_2$ in a vector space $V$ to Kronecker blocks is
defined noncanonically.  For example, let us consider
$4$-dimensional $V=\Span \{ e,p,q_1,q_2\} ,c_1=p\wedge q_1,
c_2=p\wedge q_2$. Here $V=V_1\oplus V_2$, where $V_1=\Span
\{ e\}, V_2=\Span\{ p,q_1,q_2\}$, but instead $V_1$ one can
choose any direct complement to $V_2$. However, dimensions
of the direct sums for the Kronecker blocks of equal
dimension are invariants (see~\cite{gz5},\cite{panas}). For
instance, dimension of the sum of
the trivial Kronecker blocks is equal to $\dim(\ker
c_1\cap\ker c_2)$ (see
Proposition~\ref{20.70}, below).
\end{rem}

We conclude the section by a result that will be used later
on.

\abz\label{20.70}
\begin{prop}
Let $V$ be a vector space over $\C$ and let a pair of
bivectors $c_1,c_2\in \bigwedge^2V$ be such that there are
no Jordan blocks in the decomposition of
Theorem~\ref{a20.65}. Set
\begin{eqnarray*}
\mu=\dim (\ker c_1\cap\ker c_2)\\
\mu_\la=\dim (\ker c_1\cap\ker c_\la),
\end{eqnarray*}
where $c_\la=\la_1c_1+\la_2c_2, \la_1,\la_2\not = 0$.
Then $\mu=\mu_\la$
and this number is equal to dimension of
the sum of the trivial Kronecker blocks.
\end{prop}

\noindent{\sc Proof.} Let $(V',a_1,a_2),V'\subset V,\dim
V'=2m+1, a_i\in\bigwedge^2V'$, be a nontrivial Kronecker
block.  By formula~\ref{krone}
$$
\ker
a_\la=\la_2^mq^{m+1}-
\la_1^{m-1}\la_2q^m+\cdots+(-1)^m\la_2^mq^1,
$$
where
$a_\la=\la_1a_1+\la_2a_2$ and $q^1,\ldots ,q^{m+1}$ is a
part of the basis $p^1,\ldots ,p^m,\lbr q^1,\ldots
,q^{m+1}$ in $V'^*$ dual to $p_1,\ldots p_m,\lbr
q_1,\ldots, q_{m+1}$ (let us denote these bases by ${\bf
p.},{\bf q.}$ and ${\bf p^.},{\bf q^.}$, correspondingly,
and call them adapted to the pair $a_1,a_2$).  The
above formula shows that $\ker a_1\cap\ker a_\la=\{0\}$ if
$\la_2\not =0$.

Now, let
$V=\oplus_{j=1}^kV_j,c_i=\sum_{j=1}^kc^{(j)}_i,i=1,2$, be
the decomposition to the Kronecker blocks and let
$V_{k'+1},\ldots ,V_k$ be all trivial ones.
Consider a basis of $V$ of the following form
\begin{eqnarray*}
{\bf p.}^{(1)},{\bf q.}^{(1)},
\ldots ,
{\bf p.}^{(k')},{\bf q.}^{(k')},
r_1,\ldots ,r_{k-k'},
\end{eqnarray*}
where ${\bf p.}^{(j)},{\bf q.}^{(j)}$ is a basis of $V_j$
adapted to $c_1^{(j)},c_2^{(j)},j=1,\ldots ,k'$, and
$r_1,\ldots ,r_{k-k'}$ generate $V_{k'+1},\ldots ,V_k$,
respectively. The dual basis will be of the form
\begin{eqnarray*}
{\bf p^.}^{(1)},{\bf q^.}^{(1)},
\ldots ,
{\bf p^.}^{(k')},{\bf q^.}^{(k')},
r^1,\ldots ,r^{k-k'}
\end{eqnarray*}
and the above considerations show that $\ker c_1\cap\ker
c_\la$ is generated by $r^1,\ldots ,r^{k-k'}$ if $\la_2\not
= 0$; consequently $\dim(\ker c_1\cap\ker
c_\la)$ is constant over $\la,\la_2\not = 0$ and equal to
dimension of
the sum of the trivial Kronecker blocks. q.e.d.

\section{Reductions and realizations of bihamiltonian structures.}
\label{s30}
Our next aim is to prove that a Poisson reduction
of a bihamiltonian structure is again a bihamiltonian structure.
This result follows from the naturality of  behavior of the
Schouten bracket with respect to the reduction.

\noindent\abz\label{a30.10}\begin{abza}\rm
Consider a $C^{\infty}$-smooth surjective
submersion $p:M\longrightarrow M'$
such that $p^{-1}(x')$ is connected for any $x'\in M'$.
The foliation of its leaves will be denoted by ${\cal K}$.
Write $p_*:TM\longrightarrow TM'$
for the corresponding tangent bundle morphism,
$\bigwedge^kp_*:\bigwedge^kTM\longrightarrow\bigwedge^kTM'$
%($\bigwedge^kp_*^{\C}:\bigwedge^kT^{\C}M\longrightarrow\bigwedge^kT^{\C}M'$)
 for its
exterior power extension and  $\ker \bigwedge^kp_*$
%($\ker\bigwedge^kp_*^{\C}$)
for a subbundle in $\bigwedge^kTM$
%($\bigwedge^kT^{\C}M$)
that is a kernel of $\bigwedge^kp_*$
%($\bigwedge^kp_*^{\C}$).
Multivector fields on $M$ or $M'$
will be called multivectors for short.

If $(U,\{x^1,\ldots ,x^l,y^1,\ldots,y^{m'}\})$ is a local coordinate
system on $M$ such that $m'=\dim M'$ and
$y^1,\ldots,y^{m'}$ are constant along ${\cal K}$, then the
restriction $Z|_U$ of $Z\in\Gamma(\bigwedge^kTM)$  belongs
to $\Gamma(\ker\bigwedge^kp_*)(U)$ if and only if each term
of its decomposition with respect to
$\{\d_{x^1},\ldots,\lbr\d_{x^l},\lbr\d_{y^1},\lbr\ldots,\d_{y^{m'}}\}$
contains at least
one $\d_{x^i},1\leq i\leq l$.
\end{abza}

\abz\label{a30.20}
\begin{th}
Let $Z\in\Gamma(\bigwedge^kTM)$. The following conditions are equivalent:
\begin{description}
\item[($i$)] ${\cal L}_XZ\in\Gamma(\ker\bigwedge^kp_*)\
\forall X\in\Gamma(\ker p_*)$, where ${\cal L}_X$ is a Lie derivation;
\item[($ii$)] $\phi_{t,*}^XZ-Z\in\Gamma(\ker\bigwedge^kp_*)\
\forall t\forall X\in\Gamma (\ker p_*)$, where $\phi_t^X$ denotes the flow
of the vector $X$;
\item[($iii$)] in any local coordinate system
$(U,\{x^1,\ldots ,x^l,y^1,\ldots,y^{m'}\})$
on $M$ such that $m'=\dim M'$ and
$y^1,\ldots,y^{m'}$ are constant on the leaves of $p$ the multivector $Z$
can be written as
\sabz\label{f30.19}
\begin{equation}
Z(x,y)=Z'(y)+\tilde{Z}(x,y),
\end{equation}
where
\sabz\label{f30.20}
\begin{equation}
Z'(y)= Z'^{i_1\ldots i_k}(y)\d_{y^{i_1}}\wedge\ldots\wedge\d_{y^{i_k}}
\end{equation}
and
$\tilde{Z}\in\Gamma(\ker\bigwedge^kp_*)(U)$.
\end{description}
If one of these conditions is satisfied for $Z$, then
$Z'(x')=\bigwedge^kp_*(Z(x))$, $x'\in M',x\in p^{-1}(x')$, is a correctly
defined multivector on $M'$. Moreover, if
$(p(U),\{y^1,\ldots,y^{m'}\})$ is the
induced local coordinate system on $M'$, then the corresponding
local expression for $Z'$ coincides with (\ref{f30.20}).
\end{th}

\noindent{\sc Proof.} In order to prove the last assertion it
is sufficient to note
that for any two points $x_1, x_2\in p^{-1}(x)$ there exist
$X_1,\ldots,X_s\in\Gamma(\ker p_*)$ and $t_1,\ldots,t_s\in\R$ such that
$\phi_{t_1}^{X_1}\circ\cdots\circ\phi_{t_s}^{X_s}(x_1)=x_2$ and then use the
second condition.

Obviously, $(ii)\Rightarrow (i)$. To prove the converse we choose a
vector bundle direct
decomposition $TM=\ker p_*\oplus C$ such that $Z\in \Gamma (C)$ if $Z\not
\in \Gamma
(\ker p_*)$ and $C$ is arbitrary  otherwise. Let $\Pi:\Gamma(TM)\longrightarrow
\G (C)$ be a projection on $\G (C)$ along $\G (\ker p_*)$. Then
$$\frac{d}{dt}\Pi(\phi_{t*}^XZ-Z)=\Pi\frac{d}{dt}(\phi_{t*}^XZ-Z)=
\Pi(-\phi_{t*}^X[X,Z])=0$$
(we have used the equality $\frac{d}{dt}\phi_{t*}^XZ=-\phi_{t*}^X[X,Z]$ and
the fact that $[X,Z]={\cal L}_XZ$, see \cite{lichn}).
Thus $\Pi(\phi_{t*}^XZ-Z)$ is a constant with respect to $t$ multivector and,
since $\Pi(\phi_{t*}^XZ-Z)|_{t=0}=\Pi(0)=0$, we deduce that
$\Pi(\phi_{t*}^XZ-Z)\equiv 0$.

The equivalence $(i)\Leftrightarrow (iii)$ follows from the
local expression
\sabz\label{f30.201}
\begin{equation}
[X,Z]^{i_1\ldots i_k}=\frac{1}{k!}\epsilon_{s_1\ldots
s_k}^{i_1\ldots i_k}X^r\d_rZ^{s_1\ldots s_k}-
\frac{1}{(k-1)!}\epsilon_{is_2\ldots s_k}^{i_1\ldots
i_k}Z^{rs_2\ldots s_k}\d_rX^i
\end{equation}
for the Schouten bracket (\cite{lichn}). Indeed, if one applies (\ref{f30.201})
to the local coordinate system from condition $(iii)$ one finds
that ${\cal L}_XZ\in\G(\ker \bigwedge^kp_*)$ if and only if
(\ref{f30.19}) holds. q.e.d.

\abz\label{a30.30}
\begin{defi}
We say that a multivector $Z\in\G(\bigwedge^kTM)$ is
projectable or admits the push-forward
 if one of the conditions of
Theorem \ref{a30.20} is satisfied.  The
push-forward, which will be denoted by $Z'$, is the uniquely
defined multivector from $\G(\bigwedge^kTM')$, see Theorem
\ref{a30.20}.
\end{defi}

\abz\label{a30.35}
\begin{defi} A complex multivector $Z\in\G(\bigwedge^kT^{\C} M)$
admits the push-forward $Z'\in\G(\bigwedge^kT^{\C} M')$ if the
multivectors $\Re Z,\Im Z\in\linebreak\G(\bigwedge^kTM)$ do
so. We put $Z'=(\Re Z)'+\i (\Im Z)'$.  \end{defi}

\abz\label{30.37}
\begin{cor}
Let $c$ be a (complex) bivector on $M$ admitting the
push-forward $c'\in\G (\bigwedge^2TM')$ ($c'\in\G
(\bigwedge^2T^\C M')$). Then for any $x'\in M'$ and any
$x\in p^{-1}(x')$ the following conditions hold:
\begin{description}
\item[($i$)] the subspace $p_{*,x}(á^\sharp ((T_x{\cal
K})^\bot)\subset T_{x'}M'$  ($p_{*,x}^\C (á^\sharp
((T_x^\C{\cal K})^\bot)\subset T_{x'}^\C M'$), where $\bot$
is the annihilator sign, is independent of $x$;
\item[($ii$)] the kernel of the map $p_{*,x}|_{c^\sharp
((T_{x}\K)^\bot)}$
($p_{*,x}^\C |_{c^\sharp
((T_{x}^\C\K)^\bot)}$) equals \linebreak $c^\sharp
((T_x\K)^\bot)\cap T_x\K$
($c^\sharp
((T_x^\C\K)^\bot)\cap T_x^\C\K$);
\item[($iii$)] the characteristic subspace
of the push-forward can be described by
the following isomorphism
\begin{eqnarray*} P_{c',x'}\cong
c^\sharp ((T_x\K)^\bot)/ (c^\sharp ((T_x\K)^\bot)\cap
T_x\K) \\ \mbox{\em (}P_{c',x'}\cong c^\sharp
((T_x^\C\K)^\bot)/ (c^\sharp ((T_x^\C\K)^\bot)\cap
T_x^\C\K)\mbox{\em )}.
\end{eqnarray*}
\end{description}
\end{cor}

\noindent{\sc Proof.} $iii$) follows from $i$) and $ii$).
These last are consequences of Theorem \ref{a30.20}. q.e.d.

\abz\label{30.38}
\begin{rem}\rm
Although $\dim c^\sharp ((T_x\K)^\bot)/ c^\sharp
((T_x\K)^\bot)\cap T_x\K$ is constant along $\K$,
dimensions of $c^\sharp ((T_x\K)^\bot)$ and $c^\sharp
((T_x\K)^\bot)\cap T_x\K$ may not be so. For example, let
$p:\R^4\rightarrow\R^3$ be the
projection$(x_0,x_1,x_2,x_3)\mapsto (x_1,x_2,x_3)$ and let
$c=x_0\d_{x_0}\wedge\d_{x_1}+\d_{x_2}\wedge\d_{x_3}$. Then
$c$ is projectable, but dimension of $c^\sharp
((T_x\K)^\bot)=\Span\{\d_{x_2},\d_{x_3},x_0\d_{x_0}\}$
jumps at $0$.
\end{rem}

\abz\label{a30.40}
\begin{prop}
Let a bivector $Z_i\in \G(\bigwedge^2TM)$ admit the push-forward
$Z'_i\in\G(\bigwedge^2TM'),i=1,2$. Then a trivector $Z=[Z_1,Z_2]
\in\G(\bigwedge^3TM)$ admits the push-forward $Z'\in\G(\bigwedge^3TM')$ and
$Z'=[Z'_1,Z'_2]$.
\end{prop}

\noindent{\sc Proof.} In any local coordinate system
as in condition $(iii)$ of Theorem
\ref{a30.20} $Z_i$ can be written in the form
$$Z_i(x,y)=Z'_i(y)+\tilde{Z}_i(x,y),$$
where $Z'_i(y)={Z'}_i^{jk}(y)\d_{y^j}\wedge \d_{y^k}$ and
$\tilde{Z}_i\in\G(\ker \bigwedge^2p_*)(U)$.
By formula (\ref{f10.20})
$$[Z_1,Z_2](x,y)=[Z'_1,Z'_2](y)+\tilde{Z}(x,y),$$
where $\tilde{Z}\in\G(\ker\bigwedge^3p_*)(U)$. Thus by Theorem
\ref{a30.20}  $Z$ admits
the push-forward $Z'$ and $Z'=[Z'_1,Z'_2]$. q.e.d.

\abz\label{a30.50}
\begin{cor}
Let $(c_1,c_2)$ be a Poisson pair on $M$ such that $c_i$ admits the
push-forward $c'_i\in\G(\bigwedge^2TM'),i=1,2$, and $c'_1,c'_2$ are linearly
independent.  Then $(c'_1,c'_2)$ is a Poisson
pair on $M'$.
\end{cor}

\noindent {\sc Proof}  follows immediately from
Propositions \ref{a10.190} and \ref{a30.40}. q.e.d.

\abz\label{a30.60}
\begin{cor}
Let $(M,\omega)$ be a complex symplectic manifold and let $G$ be a real Lie
group acting on $M$ by biholomorphic symplectomorphisms. Assume that $M/G$ is
a manifold. Then $c_1=\Re (\omega^{-1}),c_2=\Im (\omega^{-1})$ (see Example
\ref{a10.240}) admit the push-forwards
$c'_1,c'_2\in\G(\bigwedge^2T(M/G))$ and $(c'_1,c'_2)$ is a Poisson pair on
$M/G$, provided that $c'_1,c'_2$ are linearly independent.
\end{cor}

\noindent{\sc Proof.} It is sufficient to observe that: a)
${\cal L}_{X_j}c_i=0,i=1,2$, for generators $X_1,\ldots,X_l\in \G TM$ of
the $G$-action; b) an arbitrary vector $X\in\G(ker p_*)$, where
$p:M\longrightarrow M/G$ is a natural projection, is expressed as $X=a^jX_j$
for some $a^j\in{\cal E}(M)$ and ${\cal
L}_Xc_i=[a^jX_j,c_i]=[a^j,c_i]\wedge X_j\in \G(\ker\bigwedge^2p_*),i=1,2$ (we
have used the standard properties of the Schouten bracket, see
\cite{lichn},p.454). q.e.d.

\abz\label{a30.70}
\begin{defi}
Let $p:M\longrightarrow M'$ be as in \ref{a30.10}. A
bihamiltonian structure $(J,c_1,c_2)$ on $M$ is called
projectable (via $p$) if the bivectors $c_1,c_2$ are so and
their push-forwards $c'_1,c'_2$ are linearly independent or
zero.  The bihamiltonian structure generated by $c'_1,c'_2$
on $M'$ will be denoted by $J'$ and will be called the
push-forward or reduction of $J$.
\end{defi}

\abz\label{30.72}
\begin{defi} Let $p:M\longrightarrow M'$ be as in
\ref{a30.10} and let $J$ be a projectable bihamiltonian
structure. We say that the triple $(M,J,\K)$ is a
realization of $J'$. If, moreover, $J$ is (holomorphic)
symplectic (Definitions \ref{a10.230}, \ref{a10.247}), we
call $(M,J,\K)$ (holomorphic) symplectic realization.
\end{defi}

\section{From symplectic to complete}

Let $p:M\longrightarrow M'$ be as in
\ref{a30.10} and let $J$ be a projectable symplectic
bihamiltonian structure on $M$ with the push-forward $J'$.
In this section we discuss some conditions on the triple
$(M,J,\K)$ that guarantee the completeness of $J'$.

In view of Corollary \ref{30.37}, ($iii$) and the
definition of completeness (\ref{a20.40}) our
considerations should be linear algebraic in essence.

\abz\label{b30.10}
\begin{abza}\rm
So let $V$ be a vector space over $\C$ and let
$c_1,c_2\in\bigwedge^2V$ be such that the bihamiltonian
structure $J=\{c_\la\}_{\la\in\C^2},c_\la
=\la_1c_1+\la_2c_2, \la=(\la_1,\la_2)$, where $c_\la$ is
considered as a constant complex bivector field, is
symplectic (Definition \ref{a10.230}). Also, let $K\subset
V$ be a subspace such that the push-forwards $c'_1,c'_2\in
\bigwedge^2(V')$, where $V'=V/K$, (via the canonical
projection $p:V\rightarrow V/K$) are linearly independent.

Set $R_0=\max_{\la\in\C^2}\ra c_\la,
R'_0=\max_{\la\in\C^2}\ra c'_\la$, where
$c'_{\la}=\la_1c'_1+\la_2c'_2$, and
$$
d_\la=\dim
c_\la^\sharp(K^\bot)/(K\cap c_\la^\sharp(K^\bot)).
$$
\end{abza}

\abz\label{b30.20}
\begin{prop}
The condition of completeness
\begin{description}
\item[($**$)] $\ra c'_\la=R'_0$ for any $\la\in \C^2\smallsetminus\{(0,0)\}$
\end{description}
holds if and only if $d_\la$ is independent of $ \la\in
\C^2\smallsetminus\{(0,0)\}$.
\end{prop}

\noindent{\sc Proof.} By Corollary \ref{30.37}, ($iii$) the
space $c_\la^\sharp(K^\bot)/(K\cap c_\la^\sharp(K^\bot))$
is isomorphic to the characteristic subspace
$P_{c'_\la}=(c'_{\la})^\sharp(V')$ of the push-forward
$c'_\la$. q.e.d.

Under some additional assumption one can characterize the
condition of completeness ($**$) in terms of the subspace
$K\cap c_\la^\sharp(K^\bot)$ itself.

\abz\label{b30.30}
\begin{prop}
Let $\La_1=\Span_\C\{\hat{\la}_1\},\ldots \La_s=\Span_\C
\{\hat{\la}_s\}$ be the complex lines in $\C^2$ on which
$\ra c_\la$ is less than maximal.
Set $\La= \bigcup_{i=1}^s\La_i$ and $ k_\la= \dim K\cap
c_\la^\sharp(K^\bot)$. Assume that $K^\bot\cap \ker
c_{\hat{\la}_i}^\sharp=\{0\}, i=1,\ldots ,s$.

Then
$k_\la=\codim P_{c'_\la}$, where $P_{c'_\la}$ is the
characteristic subspace of $c'_\la$.

Consequently, the
condition ($**$) of Proposition \ref{b30.20} holds if and
only if
\begin{description}
\item[($i$)] $k_\la=k$ is
constant over $\la\in \C^2\smallsetminus\La$;
\item[($ii$)]
$k_{\hat{\la}_1}=\cdots =k_{\hat{\la}_s}=k$.
\end{description}
\end{prop}

\noindent{\sc Proof.} Since $c_\la, \la\not\in\La$ is
nondegenerate, $\codim(K+c_\la^\sharp(K^\bot))=\dim (K\cap
c_\la^\sharp(K^\bot)$. On the other hand,
$\codim(K+c_\la^\sharp(K^\bot))=\codim P_{c'_\la}$ for any
$c_\la$. Thus $k_\la=\codim P_{c'_\la}$ for
$\la\not\in\La$.

Now, the condition $K^\bot\cap \ker
c_{\hat{\la}_i}^\sharp=\{0\}$ implies the equalities
\linebreak $\dim c_{\hat{\la}_i}^\sharp(K^\bot)=\dim
c_\la^\sharp(K^\bot), i=1,\ldots ,s$, where
$\la\not\in\La$. Hence
$\codim(K+c_{\hat{\la}_i}^\sharp(K^\bot))=\dim (K\cap
c_{\hat{\la}_i}^\sharp(K^\bot)$ and
$k_{\hat{\la}_i}=\codim P_{\hat{\la}_i},i=1,\ldots ,s$.
q.e.d.

The following theorem gives the necessary and sufficient
conditions for the completeness of the reduction $J'$ of a holomorphic
symplectic bihamiltonian structure $J$ under an additional
assumption corresponding to that in
Proposition~\ref{b30.30}.  Namely, the foliation ${\cal K}$
of the leaves of the projection $p$ is supposed to be a
generic $CR$-foliation.

%It might seem surprising that $J'$ can be complete since $\ra
%c_\la, c_\la \in J$, jumps when $c_\la$ passes through the
%bivectors proportional to $c=\omega^{-1}$ and $\bar{c}$ (cf.
%Definition \ref{a20.40} and Theorem \ref{a10.245}). However,
%this can be explained by the fact that the only Casimir
%functions for $c$ ($\bar{c}$) are (anti)holomorphic functions
%and no of them are constant along ${\cal K}$.

Let $\hat{\la}_1=(1,\i), \hat{\la}_2=(1,-\i)$ and let
 $\Lambda$ denote the  cross
$\Span_\C\{\hat{\la}_1\}\cup\Span_\C\{\hat{\la}_2\}\subset \C^2$.
% $\{ (\la_1,\la_2)\in\C^2;(\la_1-(1+\i))\!\cdot\!\(la_2-(1-\i))
%=0\}\subset\C^2$.

\abz\label{ns}
\begin{th}
Let $(M,\omega)$ be a complex symplectic manifold with the
corresponding holomorphic symplectic bihamiltonian
structure $J$ (see Definition \ref{a10.247}) and let
$p:M\longrightarrow M'$ be as in \ref{a30.10}. Assume that
the foliation ${\cal K}$ is a generic $CR$-foliation on $M$ and that
$c=\omega^{-1}$
admits the push-forward $c'\in\G(\bigwedge^2T^{\C} M')$.
For $x'\in M',x\in p^{-1}(x')$, and $\la\in\C^2\smallsetminus\Lambda$ set
\begin{eqnarray*}
k_\la^{x}=\dim T_x^\C\K\cap(T_x^\C\K)^{\bot\omega_\la(x)},\\
k^x=\dim  T_x^{1,0}\K\cap(T_x^{1,0}\K)^{\bot\omega(x)},
\end{eqnarray*}
where $\omega_\la=(c_\la)^{-1}=(\la_1\Re c+\la_2\Im
c)^{-1}$. Assume that these numbers are constant along $\K$
(cf. Remark \ref{30.38}) and set $k_\la^{x'}=k_\la^x,
k^{x'}=k^x$.

Then $k_\la^{x'}=\codim_\C P_{c'_\la,x'},k^{x'}=\codim_\C
P_{c',x'}=\codim_\C
P_{\bar{c}',x'}$.

Consequently, the reduction $J'$ of $J$ via $p$ is complete
at a point $x'\in M'$ if and only if \begin{description}
\item[($i$)] $k_\la^{x'}=k$ is constant in $\la$;
\item[($ii$)] $k^{x'}=k$ \item[($iii$)] $k=\min_{y'\in
M'}k_\la^{y'}$.  \end{description}
Here $\bot\omega_\la,\bot\omega$ denote the skew-orthogonal
complements in $T^\C M,T^{1,0}M$ with respect to
$\omega_\la,\omega$, correspondingly, i.e. $(T_x^\C{\cal
K})^{\bot\omega_\la}=c_\la^\sharp((T_x^\C{\cal
K})^\bot),\linebreak (T_x^{1,0}{\cal
K})^{\bot\omega}=c^\sharp((T_x^{1,0}{\cal K})^\bot)$.
\end{th}

\noindent{\sc Proof.} If $W$ is a real vector space with
a complex structure $\J$ and $Y\subset W$ a subspace, let
$W^{1,0}$ denote the space $\{w-\i\J w; w\in W\}\subset
W^\C$ and let $Y^{1,0}=Y^\C\cap W^{1,0}=\{y-\i\J y; y\in
Y\cap\J Y\}$ (cf.~\ref{a10.145}).

Put $V=T_x^\C M, K=T_x^\C\K$. We claim
that the assumptions of Proposition~\ref{b30.30} are
satisfied.  Indeed, by Proposition~\ref{a10.245} $\La$ is
appropriate since the only, up to rescaling, degenerate
bivectors from family $J$ are $c$ and $\bar{c}$.  On the
other hand the condition $T_x\K+\J T_x\K=T_xM$ of
$CR$-genericity for $\K$ implies that $(T_x\K)^\bot\cap
\J^*((T_x\K)^\bot)=\{0\}$, where $\J^*:T^*M\rightarrow
T^*M$ is adjoint to the complex structure operator $\J$ on
$TM$. This means equalities $((T_x\K)^\bot)^{1,0}=\{
0\}=((T_x\K)^\bot)^{0,1}$   equivalent to
$K^\bot\cap T_x^{1,0}M=0=K^\bot\cap T_x^{0,1}M$. Recalling
that $T_x^{0,1}M=\ker c(x)$ and $T_x^{1,0}M=\ker\bar{c}(x)$
we get the claim.

Now, put $k_\la=k^{x'}_\la$ and
$k_{\hat{\la}_1}=k_{\hat{\la}_2}=k^{x'}$ and apply
Proposition~\ref{b30.30}. Conditions ($i$) and ($ii$) are
equivalent to the constancy of rank for $c_\la\in
J'(x'),\la\not = 0$. Its maximality is guaranteed by
($iii$). q.e.d.

\abz\label{a30.100}
\begin{cor}
In the assumptions of Theorem \ref{ns} suppose that ${\cal K}$ is
completely real (Definition \ref{a10.155}). Then if $J'$ is
nontrivial it is not complete.  \end{cor}

\noindent{\sc Proof.} Assume the contrary. By condition $(ii)$
corank of any $c'\in
J'\smallsetminus\{0\}$ is $0 $. This contradicts with the
definition of completeness. q.e.d.

Given a complete bihamiltonian structure $J'$ on $M'$,
consider all its realizations with ${\cal K}$ being a
generic $CR$-foliation. Then the smallest realizations in
this class will be characterized by the smallest
difference $T^{1,0}{\cal K}\smallsetminus T^{1,0}{\cal K}\bigcap
(T^{1,0}{\cal K})^{\bot\omega}$.

\abz\label{a30.110}
\begin{defi}
Let $J'$ be a complete bihamiltonian structure on $M'$. Its realization
$(M,\omega)$ is called minimal if $T^{1,0}{\cal K}\bigcap (T^{1,0}{\cal
K})^{\bot\omega}=T^{1,0}{\cal K}$, i.e.  ${\cal K}$ is a
$CR$-isotropic foliation (Definition \ref{a10.160}).
\end{defi}

We shall give another characterization of the minimal
realizations below.

\abz\label{a30.120}
\begin{abza}\rm
There is a natural $CR$-coisotropic foliation ${\cal L}\supset
{\cal K}$ associated with any realization $J$ on $(M,\omega)$ of
a complete $J'$. This
foliation is built as follows. Consider the "real form" $J'_{\R}$
of $J'$, i.e. the following  real bihamiltonian structure on $M'$
(cf. Proposition~\ref{a10.225})
$$
\{\la_1\Re\,c'+\la_2\Im\,c'\}_{\la=(\la_1,\la_2)\in\R^2},
$$
where $c'=p_*c,c=(\omega)^{-1}$.
Now take the bilagrangian foliation ${\cal L}'$ of $J'$
(see Definition~\ref{d20.47}). The equations for ${\cal L}'$
are the functions from the involutive family $F'_0$
(see \ref{a10.248},~\ref{a10.250}).  We define ${\cal L}$
as $p^{-1}({\cal L}')$. Note that it is $CR$-coisotropic
due to the fact that its equations $f\in p^{-1}(
F'_0)$ are in involution with respect to $c$.  \end{abza}

\abz\label{a30.130}
\begin{prop}
A realization $J$ of a complete bihamiltonian structure $J'$ is
minimal if and only if the
foliation ${\cal L}$ is a $CR$-lagrangian foliation.
\end{prop}

\noindent{\sc Proof.}  Let $2n, r$ denote rank and corank of
the bivector $c'\in J'$, respectively, and let $\dim_{\C}M=2N$.
By Definition \ref{a30.110} and Theorem \ref{ns} $J$ is
minimal if and only if $r=d$, where $d$ is $CR$-dimension of the
leaves of ${\cal K}$. On the other hand,
since
${\cal K}$ is generic, $CR$-codimension of the leaves is equal
to their real codimension, hence
$2N-d=2n+r$. Thus the minimality of $J$ is equivalent to the equality
$n+r=N$ that is necessary and sufficient for ${\cal L}$ to be
$CR$-lagrangian (see ~\ref{a10.160}).

\section{Canonical complex Poisson pair associated
with complexification of Lie algebra}

\abz\label{d40.10}
\begin{abza}\rm
Let $\A_0$ be a nonabelian  Lie
algebra over $\R$ and let $\A=\A_0 ^\C$
be its complexification. All our further results can be
formulated and proved using $\A_0,\A$ only. But we
introduce the corresponding Lie groups for the convenience.

So $G_0$ will stand for a connected simply connected Lie
group with the Lie algebra $\Lie(G_0)=\A_0$ and $G$ for
a connected simply connected complex Lie group with $\Lie
(G)=\A$. One can consider $G_0$ as a real Lie subgroup in
$G$ (see~\cite{bourb0}, III.6.10).

Write $\A_0^*, \A^*$ for the dual
spaces.
Fix a basis $e_1,\ldots,e_n$ in
$\A_0$; let $c_{ij}^k$ be the corresponding structure
constants and let $z_1=x_1+\i y_1,\ldots,z_n=x_n+\i y_n$ be
the complex linear coordinates in $\D$ associated to the
dual basis in $\D\supset\D_0$. There are the standard
linear bivectors $c=c_{ij}^kz_k\frac{\d}{\d
z_i}\wedge\frac{\d}{\d z_j}$ in $\D$ and
$c_0=c_{ij}^kz_k\frac{\d}{\d x_i}\wedge\frac{\d}{\d x_j}$
in $\A_0^*$.  They can be defined intrinsically for
instance as the maps \sabz\label{f40.10}
$$
\begin{array}{ll}
\D\stackrel{c}{\longrightarrow}\D\bigwedge\D,&
\D_0\stackrel{c}{\longrightarrow}\D_0\bigwedge\D_0
\end{array}
$$
dual to the Lie brackets
$[\,,]:\A\bigwedge\A\longrightarrow\A$ and
$[\,,]:\A_0\bigwedge\A_0\longrightarrow\A_0$.

It is well-known that the symplectic leaves of $c_0$
(respectively $c$) are the coadjoint orbits for $G_0$
(respectively $G$).  Also, there is a natural (right)
action of $G_0$ on $\D$:

$$ (a+\i b)g=Ad^*g(a)+\i Ad^*g(b), g\in
G_0, a,b\in\A_0.
$$

Let $\Sing\, \D$ be the union of
symplectic leaves of nonmaximal dimension for $c$
\end{abza}

\abz\label{d40.20}
\begin{prop}
The set $\Sing\, \D$ is algebraic.
%Moreover, $\Sing\,\D$ is the complexification of
%$\Sing\,\D_0$, i.e. the defining polynomials for
%$\Sing\,\D$ have real coefficients and are obtained
%from the defining polynomials of $\Sing\,\D_0$ by
%substitution $\{x_1,\ldots ,x_n\}\mapsto\{z_1,\ldots
%,z_n\}$. Consequently, $\codim_\C\Sing\,\D\le
%\linebreak
%\codim_\R\Sing\,\D_0$.
\end{prop}

\noindent{\sc Proof.} The defining  polynomials for
$\Sing\,\D$ are minors of $m$-th order of the
$n\times n$-matrix $||c_{ij}^kz_k||$, where $m=\ra c$.
q.e.d.

\abz\label{d40.15}
\begin{conv}\rm In the sequel we
shall assume that the nonabelian Lie algebra $\A$
satisfies condition $\codim_\C\Sing\,\D\ge 3$.
\end{conv}

\abz\label{d40.30}
\begin{abza}\rm
This condition is satisfied by a wide class of Lie algebras
including the semisimple ones. Indeed, in the semisimple
case we can identify $\D$ and $\A$ by means of the Killing
form. On the other hand, it is well known that  the
algebraic set of all nonregular (regular means semisimple
contained in the unique Cartan subalgebra) elements is at
least of codimension three and contains $\Sing\,\D$.
\end{abza}

\abz\label{d40.150}
\begin{defi}
Let us introduce a set
$$
{\cal C}=\{z\in\D;\exists
(\la_1,\la_2)\in\C^2\smallsetminus\{(0,0)\},\la_1z+\la_2\bar{z}\in \Sing\,\D\},
$$
where the bar stands for the complex conjugation
corresponding to $\A_0\subset\A$, and call it the
incompleteness set (see~\ref{a40.100} for the explanation
of  this terminology).
\end{defi}

\abz\label{a40.20}
\begin{prop}
 The incompleteness set ${\cal
C}$ is a real algebraic set of positive codimension.
\end{prop}

\noindent{\sc Proof.} We use the product $\Pi=\D\times
(\C^2\smallsetminus\{(0,0)\})$ with
the coordinates $z_1,\ldots,z_n,\la_1,\la_2$ and the real
algebraic map $\phi:\Pi\longrightarrow \D$ given by the
formula $$ (z_1,\ldots,z_n,\la_1,\la_2)\mapsto
(\la_1z_1+\la_2\bar{z_1},\ldots,\la_1z_n+\la_2\bar{z_n}).
$$
The set ${\cal C}$ can be regarded as $pr_1(\phi^{-1}(\Sing\,\D))$, where
$pr_1$ is the projection onto $\D$.

The above construction shows that $\dim_\R{\cal
C}\le\dim_\R\Sing\,\D+4$ q.e.d.

\abz\label{d40.200}
\begin{exa}\rm
Let $\A_0=so(3), \A=sl(2,\C)\cong\C^3$. Then
$\Sing\,\D=\{0\}$, ${\cal C}=\{z\in\D;z\, \mbox{linearly
independent with}\, \bar{z}\}$; consequently ${\cal C}$ is
described by two real equations:
$z_1\bar{z}_2-z_2\bar{z}_1=0,z_1\bar{z}_3-z_3\bar{z}_1=0$.
The set $\{z\in\D;z\, \mbox{linearly
independent with}\, \bar{z}\}$ is contained in ${\cal C}$
for arbitrary $\A$.  \end{exa}

\abz\label{d40.21}
\begin{abza}\rm
Now, we shall introduce a remarkable pair of complex
bivectors on $\D$ playing the crucial role in the sequel of
the paper.  This pair is $(c,\tilde{c})$, where $c$ is as
in \ref{d40.10} and $\tilde{c}$ is given by
$\tilde{c}=c_{ij}^k\bar{z}_k\frac{\d}{\d z_i}\wedge\frac{\d}{\d
z_j}$. One can define $\tilde{c}$ intrinsically by the diagram
$$
\begin{array}{lcc}
\D & \stackrel{c}{\longrightarrow} & \D\bigwedge\D\\
\uparrow \bar{\cdot} & & \uparrow \tilde{c}\\
\D & = & \D,
\end{array}
$$
where $c$ is from (\ref{f40.10}) and $\bar{\cdot}$ stands for the complex
conjugation corresponding to the real form $\A_0\subset\A$.
\end{abza}

\abz\label{a40.100}
\begin{prop}
\begin{description}
\item[($i$)] $\tilde{c}$ is $G_0$-invariant;
\item[($ii$)] $(c,\tilde{c})$ is a complex Poisson pair;
\item[($iii$)] $(c,\tilde{c})$ is complete  at any point $z\in \D\smallsetminus
{\cal C}$ (see Definition \ref{a20.40}).
\end{description} \end{prop}

\noindent{\sc Proof.} ($i$) follows from the $G$-invariance
of the bivector $c$ and $G_0$-equivariance of
$\bar{\cdot}$. ($ii$) is obtained by direct
calculations. The last assertion follows from Proposition
\ref{a40.20} since the set ${\cal C}$ consists precisely of
the points of incompleteness for $(c,\tilde{c})$.
Indeed, $\ra (\la_1c+\la_2\tilde{c})(z)=\ra
c_{ij}^k(\la_1z_k+\la_2\bar{z}_k)\frac{\d}{\d
z_i}\wedge\frac{\d}{\d z_j}$ is less than maximal if and
only if $\la_1z+\la_2\bar{z}\in\Sing\,\D$. q.e.d.

\abz\label{can}
\begin{defi}
The bihamiltonian structure generated by $c,\tilde{c}$ will
be denoted by $\tilde{J}$ and will be called the canonical
bihamiltonian structure.
\end{defi}

The end of this section is devoted to the study of the
first integrals (Definition~\ref{20.05}) for the canonical
bihamiltonian structure $(\tilde{J},c,\tilde{c})$.

\abz\label{d40.430}
\begin{defi}
Let $r=\cora c$ (codimension of symplectic leaf of maximal
dimension). Let us write $\ra\,\A$ for $r$ and call this
number the rank of $\A$.
\end{defi}

Note that for the semisimple case this notion of rank
coincides with the standard one, i.e. with dimension of a
Cartan subalgebra.

\abz\label{d40.43}
\begin{defi} Let
$Z_c^{hol}(U)$ denote the space of holomorphic Casimir
functions for $c$ over an open set $U\subset\D$.

An open set $U\subset\D\smallsetminus\Sing\,\D$ is called
admissible if there exist $r=\ra\,\A$ functionally
independent functions from $Z_c^{hol}(U)$.
\end{defi}

\abz\label{d40.45}
\begin{prop}
Let a set $U$ be admissible.
Given a function $g\in Z_c^{hol}(U)$, define a function
$\tilde{g}\in \E^\C(U)$ by the formula
$\tilde{g}=\overline{(\frac{\partial g}{\partial z_i})}z_i$.
Then the space
$Z_{\tilde{c}}(U)$ of (smooth) Casimir functions for
$\tilde{c}$ is equal to $\{\tilde{g};g\in
Z_c^{hol}(U)\}\cup \overline{\cal O}(U)$, where
$\overline{\cal O}(U)$ is the space of antiholomorphic
functions over $U$.
\end{prop}

\noindent{\sc Proof.} The following calculation shows that
$\tilde{g}\in Z_{\tilde{c}}(U)$:
$$
\tilde{c}(\tilde{g})_j=c_{ij}^k\bar{z}_k\frac{\partial
\tilde{g}}{\partial
z_i}=c_{ij}^k\bar{z}_k\overline{(\frac{\partial g}{\partial
z_i})}=\overline{c(g)}_j=0
$$
(here $v_j$ stands for the $j$-th component of a vector
field $v=v_i\frac{\partial }{\partial z_i}$).

Now, let $g_1,\ldots ,g_r\in Z_c^{hol}(U)$ be functionally
independent.
We note that the $(1,0)$-differentials $\d
\tilde{g}_1,\ldots ,\d\tilde{g}_r$ are linearly independent
precisely at those points where $\d g_1,\ldots ,\d g_r$
are. Thus by the dimension arguments ($\ra \tilde{c}=\ra
c$) the functions $\tilde{g}_1,\ldots ,\tilde{g}_r$
together with the antiholomorphic functions functionally
generate the space $Z_{\tilde{c}}(U)$. q.e.d.

\abz\label{d40.46}
\begin{defi}
Define $\phi_\la:\D\rightarrow \D,\la=(\la_1,\la_2)\in
\C^2\smallsetminus\{(0,0)\}$ by $\phi_\la(z)=\la_1
z+\la_2\bar{z}$. This is  a $\R$-linear isomorphism if
$|\la_1|\not = |\la_2|$ and an epimorphism onto an
$n$-dimensional ($n=\dim_\C\A$) real subspace otherwise.

An open set $U\subset\D\smallsetminus {\cal C}$ is called
$\la$-admissible if the set $\phi_\la(U)$ has an
admissible neighbourhood.

An open set $U\subset\D\smallsetminus {\cal C}$ is called
strongly admissible if it is $\la$-admissible for any
$\la\in\C^2\smallsetminus\{(0,0)\}$.
\end{defi}

\abz\label{d40.50}
\begin{prop}
Let a set $U\subset\D$ be $\la$-admissible and let $U_\la$
be an admissible neighbourhood of $\phi_\la(U)$.

Then  the space $Z_{c_\la}(U)$ of (smooth)
Casimir functions for
$c_\la=\la_1c+\la_2\tilde{c},\la_1\not = 0$, is equal to
$\{g\circ\phi_\la|_U;g\in Z_c^{hol}(U_\la)\}\cup
\overline{\cal O}(U)$.
\end{prop}

\noindent{\sc Proof.} Again, let $g_1,\ldots ,g_r\in
Z_c^{hol}(U_\la)$ be functionally independent. Obviously,
the functions $g_{\la,1}=g_1\circ\phi_\la,\ldots
,g_{\la,r}=g_r\circ\phi_\la$ are Casimir functions for
$c_\la$.  They are functionally independent on $U$ since
the Jacobi matrices $D=\frac{\d (g_1,\ldots,g_r)}{\d
(z_1,\ldots,z_n)}$ and $D_\la=\frac{\d
(g_{\la,1},\ldots,g_{\la,r})}{\d (z_1,\ldots,z_n)}$ are
related as follows
$$
D_\la(z)=\la_1D\circ\phi_\la (z).
$$
So $g_{\la,1},\ldots ,g_{\la,r}$ and $\overline{\cal O}(U)$
generate $Z_{c_\la}(U)$. q.e.d.

The following proposition shows that strongly admissible
sets exist and describes all of them in the semisimple case.

\abz\label{d40.55}
\begin{prop}
\begin{description}
\item[($i$)] Let $||  \cdot ||  $ be a norm in $\D$ Then any
open set $U\subset\D\smallsetminus\Sing\,\D$ with a sufficiently
small diameter $\diam U=\sup_{z,z'\in U}||  z-z'||  $ is
strongly admissible.  \item[($ii$)] Assume that $\A$ is
semisimple.  Then any open set $U\subset\D\smallsetminus {\cal
C}$ is strongly admissible.  \end{description} \end{prop}

\noindent{\sc Proof.} ($i$) We start from the following
claim: if a set $U$ is admissible, then the set $\la
U=\{\la u;u\in U\}$ is so for any $\la\in \C\smallsetminus\{
0\}$. Indeed, the bivector $c$ is homogeneous with
homogeneity degree
$1$: $h_{\la,*}c=\la c$, where $h_\la(z)=\la z$. Hence, if
$g_1,\ldots ,g_r$ are independent Casimir functions for $c$
over $U$, then $(h_\la^*)^{-1}g_1,\ldots ,(h_\la^*)^{-1}g_n$
are so over $h_\la(U)=\la U$.

Now, assume that the norm is so chosen that
$||  z||  =||  \bar{z}||  $. Then the inequality
$$
||  \la_1z+\la_2\bar{z}-\la_1z'-\la_2\bar{z}'||  \le |\la_1|\,
||  z-z'||  +|\la_2|\, ||  \bar{z}-\bar{z}'||
$$
shows that
\sabz\label{fdiam}
\begin{equation}
\diam\phi_\la(U)\le(|\la_1|+|\la_2|)\diam U,
\end{equation}
where $\phi_\la$ is from Definition \ref{d40.46}.

Next, choose a point $z\in U$ (note that $z$ is linearly
independent with $\bar{z}$, see Example~\ref{d40.200}) and
consider the map $\C^2\ni\la\mapsto\phi_\la(z)\in\D$. The
image of the unit sphere $S^1=\{|\la_1|^2+|\la_2|^2=1\}$
under this map can be covered by a finite number of
admissible balls $B_1,\ldots ,B_m$. Inequality \ref{fdiam}
shows that shrinking $U$ if needed one can get the
following $$ \phi_\la
(U)\subset\cup_{i=1}^mB_i\,\forall\la\in S^1.
$$
Hence, for
sufficiently small $U\ni z$ the set $\phi_\la (U)$, where
$\la\in S^1$, possesses an admissible neighbourhood and by
the above proved claim the same is true for $\phi_\la
(U),\la\not = 0$. Since all norms on $\D$ are equivalent
this completes the proof.

\noindent ($ii$) It is enough to note that there exists
a set $g_1(z),\ldots,g_r(z),r=\ra\A$ of global
holomorphic Casimir functions for $c$ that are functionally
independent on $\D\smallsetminus\Sing\, \D$.
One can identify $\A$ and $\D$ by means of the Killing form
and take for  $g_1,\ldots ,g_r$ an algebraic basis  of the
ring of $G$-invariant polynomials on $\A$. The functional
independence of these functions on $\D\smallsetminus\Sing\, \D$
is established in Theorem 0.1 of ~\cite{kost}. q.e.d.

We summarize the above results in the following
Proposition.

\abz\label{d40.57}
\begin{prop}
Let $U\subset \D$ be a strongly admissible set and let
$U_\la$ be an admissible neighbourhood of $\phi_\la (U)$.
The set of first integrals (see Definition~\ref{20.05})
$F_0(U)$ of $(\tilde{J},c,\tilde{c})$ over $U$ is
(functionally) generated by the functions from the sets
$F_1(U)$ and $\overline{\cal O}(U)$, where the last one is
the set of antiholomorphic functions on $U$ and $F_1(U)$ is
in turn generated by $Z_c^{hol}(U),\{\tilde{g};g\in
Z_c^{hol}(U)\},\{g\circ\phi_\la|_U;g\in Z_c^{hol}(U_\la)\},
\la=(\la_1,\la_2)\in\C^2,\la_1,\la_2\not =0$.
\end{prop}

\noindent{\sc Proof} follows from
Propositions~\ref{d40.45},\ref{d40.50} and from the
definition of the set $F_0(U)$. q.e.d.

The following proposition will be crucial in the proof of
our main result (Theorem~\ref{main}). As usual, given a
subspace $V\subset (T_z^\C\D)^*$, we set $V^{1,0}=V\cap
(T_z^{1,0}\D)^*$.

\abz\label{d40.60}
\begin{prop}
Let
$$
\mu(z)=\dim (\ker c(z))^{1,0}\cap(\ker\tilde{c}(z))^{1,0}
$$
and let
$$
\mu_\la(z)=\dim (\ker c(z))^{1,0}\cap(\ker c_\la(z))^{1,0}
$$
for
$c_\la=\la_1c+\la_2\tilde{c},\la=(\la_1,\la_2)\in\C^2,
\la_1, \la_2\not = 0$.

Then
\begin{description}
\item[($i$)] $\mu (z)=\mu_{\la}(z)$ for any $z\in
\D\smallsetminus{\cal C}$;
\item[($ii$)] there exists a real algebraic
set ${\cal R},{\cal C}\subset{\cal R}\subset\D$, where ${\cal C}$
is the incompleteness  set (see Definition~\ref{d40.150}),
such that $\mu(z)=\mu$ is constant and minimal on
$\D\smallsetminus{\cal R}$ and any set with these properties
contains ${\cal R}$;
\item[($iii$)] if $\A$ is semisimple
the set ${\cal R}\smallsetminus{\cal C}$ is empty and $\mu
(z)\equiv 0$ on $\D\smallsetminus{\cal C}$.
\end{description}
\end{prop}

\noindent{\sc Proof.} ($i$)
We shall use the completeness of the
bihamiltonian structure
$(J,c,\tilde{c})$ at any $z\in\D\smallsetminus{\cal C}$ (see
proof of Proposition~\ref{a40.100}). By
Theorem~\ref{a20.65} and Corollary~\ref{a20.67} the pair
$c(z),\tilde{c}(z)\in \bigwedge^2T_z^{1,0}\D$ does not have
the Jordan blocks in its decomposition. Thus we can use
Proposition~\ref{20.70} to deduce
that $\mu(z)=\mu_{\la}(z)$.

\noindent ($ii$) To prove this condition it is sufficient
to note that the subspace $(\ker c(z))^{1,0}\cap(\ker
\tilde{c}(z))^{1,0}$ annihilates the sum of characteristic
subspaces $P_{c,z}+P_{\tilde{c},z}$. Put ${\cal
R}=\{z\in\D; \dim(P_{c,z}+P_{\tilde{c},z})<m\}$, where
$m=\max_z\dim (P_{c,z}+P_{\tilde{c},z})$.  The
defining polynomials for ${\cal R}$ are the minors of
$m$-th order of the $2n\times n$-matrix
$$
\left|\left|
\begin{array}{c}
c_{ij}^kz_k \\ c_{ij}^k\bar{z}_k
\end{array}\right|\right| .
$$

If the  set ${\cal R}$ defined above lies in ${\cal C}$,
let us change the definition and put ${\cal R}={\cal C}$.

It remains to prove the inclusion ${\cal C}\subset{\cal
R}$ in the case ${\cal C}\not \supset {\cal R},{\cal C}\not
= {\cal R}$.  Introduce a set ${\cal R}_\la=\{z\in\D; \dim
(P_{c,z}+P_{c_\la,z})<m_\la\}$, where $m_\la=\max_z\dim
(P_{c,z}+P_{c_\la,z}),\la_2\not = 0$.  Then by ($i$) ${\cal
R}_\la,\la_2\not = 0$ coincides with ${\cal R}$ outside
${\cal C}$. Since ${\cal R}=Cl({\cal R}\smallsetminus{\cal C})$
and  ${\cal R}_\la=Cl({\cal R}_\la\smallsetminus{\cal C})$
(Zarisski closures), one gets ${\cal R} ={\cal R}_\la$.

Now, let $z\in{\cal C}$ and let
$\la=(\la_1,\la_2)\in\C^2,\la_2\not = 0$ be such that $\ra
c_\la(z)<R_0$ (cf.  the definition of
completeness, \ref{a20.40}). Then \linebreak
$\dim(P_{c,z}+P_{c_\la,z})<m_\la$. Consequently, ${\cal
C}\subset{\cal R}_\la={\cal R}$.
If the only (up to the
proportionality) bivector of nonmaximal rank in the family
$\{c_\la (z)\}$ is $c$, then $\dim
(P_{c,z}+P_{\tilde{c},z})<m$ and again ${\cal
C}\subset{\cal R}$.

\noindent ($iii$) It follows from
Proposition~\ref{20.70} that $\dim(\ker
c(z))^{1,0}\cap(\ker\tilde{c}(z))^{1,0}$ is equal to
dimension of the sum of the trivial Kronecker blocks for
the pair $c(z),\tilde{c}(z)$. If we consider $c,\tilde{c}$
as elements of $\G(\bigwedge^2T^{1,0}\D)$,  the same
arguments as for the method of argument translation
(\cite{panas}, Theorem 4.1) show that in the semisimple
case the trivial Kronecker blocks are absent for any $z\in
\D\smallsetminus {\cal C}$.  q.e.d.

\abz\label{d40}
\begin{defi}
We call the set ${\cal R}$  from
Proposition~\ref{d40.60} the Kronecker irregularity set and
the number $\mu$ the trivial Kronecker dimension of the Lie
algebra $\A$.
\end{defi}

The following example shows that for nonsemisimple Lie
algebras the set ${\cal R}\smallsetminus{\cal C}$ may be
nonempty and the trivial Kronecker dimension may be
nonzero.

\abz\label{d45}
\begin{exa}\rm
Let $\A=\Span\{p_1,p_2,q_1,\ldots ,q_4,f_1,\ldots
,f_4,g_1,\ldots ,g_4\}$ be a fourteen-dimensional Lie
algebra with the standard linear Poisson bivector
$c=\frac{\partial }{\partial p_1}\wedge(f_1\frac{\partial
}{\partial q_1}+\cdots + f_4\frac{\partial }{\partial
q_4})+\frac{\partial }{\partial
p_2}\wedge(g_1\frac{\partial }{\partial q_1} +\cdots +
g_4\frac{\partial }{\partial q_4})$. Then ${\cal R}$ is
given by one real equation
$$\left|
\begin{array}{cccc}
f_1 & f_2 & f_3 & f_4 \\
g_1 & g_2 & g_3 & g_4 \\
\bar{f}_1 & \bar{f}_2 & \bar{f}_3 & \bar{f}_4 \\
\bar{g}_1 & \bar{g}_2 & \bar{g}_3 & \bar{g}_4
\end{array}\right| =0.
$$

The set $\Sing\,\D$ consists of the points where the vectors
$f_1\frac{\partial
}{\partial q_1}+\cdots + f_4\frac{\partial }{\partial
q_4}, g_1\frac{\partial }{\partial q_1} +\cdots +
g_4\frac{\partial }{\partial q_4}$ are linearly dependent,
i.e. the defining equations for $\Sing\,\D$ are
$f_1g_2-f_2g_1=0, f_1g_3-f_3g_1=0,f_1g_4-f_4g_1=0$.
However, the proof of
Proposition~\ref{a40.20} shows that $\codim_\R{\cal
C}\ge\codim_\R\Sing\,\D-4$; consequently, in our example
$\codim_\R{\cal C}\ge 6-4=2$ and ${\cal C}\not ={\cal R}$.

Here
$\mu=8$ since $f_1,\ldots ,f_4,g_1,\ldots ,g_4$ are the
common Casimir functions for $c\tilde{c}$.

Also, $\mu$ will be nonzero for all
reductive nonsemisimple Lie algebras.

Note that the above examples agree with our
Convention~\ref{d40.15}.
\end{exa}

\section{$CR$-geometry of real coadjoint orbits}

We retain the notations and conventions from the previous
section. The reader is referred to Section 1 for the
$CR$-geometric concepts used below.

\abz\label{d50}
\begin{prop}
\begin{description}
\item[($i$)] The bivectors $c_1=\Re c, c_2=\Im c$
are Poisson.
\item[($ii$)] The coadjoint action of $G_\R$, where $G_\R$
is $G$ considered as a real Lie group, is hamiltonian with
respect to $c_1,c_2$ (see Definition~\ref{10.17}).
\item[($iii$)]
The generalized distribution of subspaces tangent to the
$G_0$-orbits is generated by the vector fields
$c(z_i)+\overline{c(z_i)},i=1,\ldots ,n$.
\end{description}
\end{prop}

\noindent{\sc Proof.} ($i$) The more general statement
that the pair $c_1,c_2$ is Poisson is proved by the same
arguments as in Example~\ref{a10.240}.

\noindent ($ii$) First, we shall prove that the holomorphic
coadjoint action of $G$ on $\D$ is hamiltonian in
holomorphic sense with respect to $c$. Consider the
antirepresentation $Ad^*:G\rightarrow \D$. The
corresponding Lie algebra action $\rho:\A\rightarrow
\Vect^{hol}(\D)$, where $\Vect^{hol}(\D)$ is the Lie
algebra of holomorphic vector fields on $\D$, can be
described as follows. The vector field $\rho (v),v\in\A$, is
equal to
$$
z\mapsto ad^*(v)z:\D\rightarrow \D\cong
T^{1,0}_z\D.
$$
On the other hand, if $e_1,\ldots ,e_n, c_{ij}^k$ are as
in~\ref{d40.10}, then
$$
<ad^*(e_i)z,e_j>=<z,ad(e_i)e_j>=\lbr
<z,[e_i,e_j]>=<z,c_{ij}^ke_k>=c_{ij}^kz_k,
$$
where we put $z=z_1e^1+\cdots +z_ne^n$. Hence, $\rho
(e_i)=c_{ij}^kz_k\frac{\partial }{\partial z_j}=c(z_i)$
and the corresponding antihomomorphism $\psi:\A\rightarrow
{\cal O}(\D)$ (cf. Definition~\ref{10.17}) is defined by
$e_i\mapsto z_i, i=1,\ldots ,n$.

Now, the Lie algebra action
$\rho_\R:\A_\R\cong\A_0\oplus\i\A_0 \rightarrow\Vect (\D)$,
where $\cong$ is over $\R$, corresponding to the
antirepresentation $Ad^*:G_\R\rightarrow \D$ is described
by the formulas \begin{eqnarray*}
<ad^*(e_i)z,e_j>=c_{ij}^kx_k,
<ad^*(\i e_i)z,e_j>=c_{ij}^ky_k\\
<ad^*(e_i)z,\i e_j>=c_{ij}^ky_k,
<ad^*(\i e_i)z,\i e_j>=-c_{ij}^kx_k.
\end{eqnarray*}
Consequently,
\begin{eqnarray*}
\rho_\R(e_i)=c_{ij}^kx_k\frac{\partial
}{\partial x_j}+c_{ij}^ky_k\frac{\partial
}{\partial
y_j}=
(c+\bar{c})(z_i+\bar{z_i})=(c-\bar{c})(z_i-\bar{z_i})
\end{eqnarray*}
and
\begin{eqnarray*}
\rho_\R(\i e_i)=c_{ij}^ky_k\frac{\partial }{\partial
x_j}-c_{ij}^kx_k\frac{\partial }{\partial
y_j}=
(c+\bar{c})(z_i-\bar{z_i})=(c-\bar{c})(z_i+\bar{z_i}).
\end{eqnarray*}

\noindent ($iii$) This condition follows from the proof
of ($ii$) and from the obvious equality
$(c+\bar{c})(z_i+\bar{z_i})=c(z_i)+\overline{c(z_i)}$.
q.e.d.

\abz\label{a40.50}
\begin{prop}
Let $\O$ be a $G_0$-orbit through $z_0\in \D$.
Then
$\O$ is a generic $CR$-manifold in the $G$-orbit $G(z_0)$;
\end{prop}

\noindent{\sc Proof.} The $G_0$-invariance of the complex structure ${\cal J}$ on $\D$
implies the constancy of $\dim T_z\O\bigcap{\cal J}T_z\O, z\in\O$. To prove
the genericity we note that the tangent bundle $T\O$ is generated by the
vector fields $c(z_j)+\overline{c(z_j)},j=1,\ldots,n$
(Proposition~\ref{d50},($iii$)), and that ${\cal J}$ acts
on them as follows \sabz\label{f40.50}
\begin{equation}
{\cal
J}(c(z_j)+\overline{c(z_j)})=\i(c(z_j)-\overline{c(z_j)}).
\end{equation}
Thus $T\O+{\cal J} T\O$ is generated by the real and imaginary parts of the
vector fields $c(z_j),j=1,\ldots,n$, spanning $T^{1,0}G(z_0)$. Hence
$T_z\O+{\cal J}T_z\O=T_z^{\C}G(z_0)$. q.e.d.

The next proposition gives some
characterization (another one can be found in~\ref{CR}) of
$CR$-dimension of a $G_0$-orbit.

\abz\label{40.51}
\begin{prop}
Let $\O$ be a $G_0$-orbit through $z_0\in \D$.
Write $G^z$ (respectively $G_0^z$) for the
stabilizer of $z\in\D$ in $G$ (respectively $G_0$) and
$\A^z$ ($\A_0^z$) for the corresponding Lie algebra.  Then
for any $z\in{\cal O}$
$$
T_z^{1,0}\O
\cong\A^z/(\A_0^z)^{\C}.
$$
\end{prop}

\noindent{\sc Proof.} We study the embedding of $T_z\O$ in
the tangent space $T_zG(z_0)$ to
a $G$-orbit $G(z_0)$.  At the Lie algebra level it is equal
to a map
$$
\iota:\A_0/\A_0^z\longrightarrow \A/\A^z
$$
induced by the inclusions
$\A_0\hookrightarrow\A,\A_0^z\hookrightarrow\A^z$.
The intersection $\iota(\A_0/\A_0^z)\cap{\cal
J}\iota(\A_0/\A_0^z)$, where ${\cal J}$ is the complex
structure on $\A/\A^z$ induced by the multiplication by
$\i$, is equal to
$$
\{s+\A_0^z;s\in\A_0,\exists t\in\A_0,s-\i
t\in\A^z\}.
$$

We
define a map
$$\phi:\A^z\longrightarrow\iota(\A_0/\A_0^z)\cap{\cal
J}\iota(\A_0/\A_0^z) \cong T_z\O\bigcap{\cal J}T_z\O
$$
by the formula
$$
\A^z\ni s-\i t\mapsto s+\A^z,
$$
where $s,t\in \A_0$. The kernel of $\phi$ is equal to $\{s-\i
t\in\A^z;s\in\A^z\}=\{s-\i
t\in\A^z;s\in\A_0^z=\A^z\cap\A_0\}$. Since
$ad^*s(x)=ad^*s(y)=0$ for $s\in\A_0^z,x=\Re\,z,y=\Im\,z$, one
gets
\begin{eqnarray*}
0=ad^*(s-\i t)(x+\i y)=-\i ad^*(t)(x+\i y)\Rightarrow \\
ad^*(t)x=ad^*(t)y=0\Rightarrow  t\in\A_0^z.
\end{eqnarray*}
Thus $\ker\phi=(\A_0^z)^{\C}$. The surjectivity of $\phi$
is obvious. q.e.d.

\abz\label{a40.70}
\begin{prop}
Let $\O$ be a $G_0$-orbit through
$z_0\in\D\smallsetminus\Sing\,\D$ and let $g_1,\ldots
,g_r\in Z_c^{hol}(U),r=\ra \A$, be independent
holomorphic Casimir functions for  $c$ in some
neighbourhood $U\subset\D\smallsetminus\Sing\,\D$ of $z_0$. Then
$T^{1,0}\O$ is generated over $U$ by the vector fields
$c(\tilde{g_1}),\ldots,c(\tilde{g_r})$ (see
Proposition~\ref{d40.45} for the notations).
\end{prop}

\noindent{\sc Proof.} By formula (\ref{f40.50})
$T\O\bigcap{\cal J}T\O$ is generated by linear combinations
$\alpha^j(c(z_j)+\overline{c(z_j)}), \alpha^j\in\E(\D)$,
such that there exist $\beta^j\in\E(\D)$ satisfying the
equality $$ \alpha^j(c(z_j)+\overline{c(z_j)})=\beta^j\i
(c(z_j)-\overline{c(z_j)}).  $$ This implies $$
(\alpha^j+\i\beta^j)c(z_j)+(\alpha^j-\i\beta^j)\overline{c(z_j)}=0
$$
and
\sabz\label{f40.70}
\begin{equation}
\gamma^jc(z_j)=0,
\end{equation}
where we put $\gamma^j=\alpha^j+\i\beta^j$.

In order to calculate all vector functions
$\gamma=(\gamma^1,\ldots,\gamma^n)$ satisfying (\ref{f40.70})
one observes two facts. First, that the vector functions
$\gamma_{(m)}=(\frac{\d g_m}{\d z_1},\ldots,\frac{\d
g_m}{\d z_n}), m=1,\ldots,r$, satisfy (\ref{f40.70}).
Second, the dimension arguments show that any $\gamma(z)$
for which (\ref{f40.70}) holds is a linear combination of
$\gamma_{(1)}(z),\ldots,\gamma_{(r)}(z)$ if
$z\in U$.

In other words, $T^{CR}\O$ is generated by $(\frac{\d g_m}{\d
z_j}+\frac{\d \bar{g}_m}{\d
\bar{z}_j})(c(z_j)+\overline{c(z_j)}),i=1,\ldots,r$, and $T^{1,0}\O$ by
\begin{eqnarray*}
(\frac{\d g_m}{\d
z_j}+\overline{(\frac{\d g_m}{\d
z_j})})(c(z_j)+\overline{c(z_j)})-\i{\cal J}(\frac{\d g_m}{\d
z_j}+\overline{(\frac{\d g_m}{\d
z_j})})(c(z_j)+\overline{c(z_j)})=\\
(\frac{\d g_m}{\d
z_j}+\overline{(\frac{\d g_m}{\d
z_j})})(c(z_j)+\overline{c(z_j)}+c(z_j)-\overline{c(z_j)})=\\
2c(g_m)+2\overline{(\frac{\d g_m}{\d z_j})}c(z_j)=2c(\tilde{g}_m).
\end{eqnarray*}
q.e.d.

In the next result we describe dimension and once more
$CR$-dimension of generic $G_0$-orbits.

\abz\label{CR}
\begin{cor}
Let ${\cal O}$ be any $G_0$-orbit lying in the complement
to the Kronecker irregularity set ${\cal R}$
(see Definition~\ref{d40}), which is $G_0$-invariant by
Proposition~\ref{a40.100}.  Then
\begin{description}
\item[($i$)] $\dim_\C T_z^{1,0}{\cal O},
z\in{\cal O}$, equals $r-\mu$, where $r=\ra \A$ and $\mu$
is the trivial Kronecker dimension; in particular, if $\A$
is semisimple, $\dim_\C  T_z^{1,0}{\cal O}=r$;
\item[($ii$)] $\dim {\cal O}=n-\mu$, where $n=\dim\A$.
\end{description} \end{cor}

\noindent{\sc Proof.} ($i$) follows immediately from
Propositions~\ref{a40.70}, \ref{d40.45}, and \ref{d40.60}.
($ii$) is a consequence of ($i$) and
Proposition~\ref{40.51} (since $\dim_\C\A^z=r$).  q.e.d.

The following corollary characterize generic $G_0$-orbits
in $\D$ from the symplectic point of view.

\abz\label{a40.80}
\begin{cor}
Let $\O$ be a $G_0$-orbit through
$z_0\in\D\smallsetminus\Sing\,\D$. Then $\O$ is a
$CR$-isotropic submanifold in $M$ (Definition
\ref{a10.160}).
\end{cor}

\noindent{\sc Proof.} First we shall show the
$G_0$-invariance of the functions
$\tilde{g_1},\ldots,\tilde{g_r}$. We use the equality
\sabz\label{f40.80}
\begin{equation}
c_{ij}^l\frac{\d g_m}{\d z_i}+ c_{ij}^kz_k\frac{\d^2g_m}{\d
z_i\d z_l}=0
\end{equation}
obtained by the
differentiation of the equality $c_{ij}^kz_k\frac{\d g_m}{\d
z_i}=0$ with respect to  $z_l$.  Conjugating
(\ref{f40.80}) and multiplying by $z_l$ one gets
\begin{eqnarray*}
0=c_{ij}^lz_l\overline{(\frac{\d g_m}{\d z_i})}+
c_{ij}^k\bar{z}_k\overline{(\frac{\d^2g_m}{\d z_i\d z_l})}z_l=\\
(c_{ij}^kz_k\frac{\d}{\d z_i}+c_{ij}^k\bar{z}_k\frac{\d}{\d
\bar{z}_i})[\overline{(\frac{\d g_m}{\d z_l})}z_l]=\\
-(c(z_j)+\overline{c(z_j)})\tilde{g}_m.
\end{eqnarray*}
This proves the claim.

Now, recall that $(T^{1,0}\O)^{\bot\omega}$ is generated by
the vector fields $c(f)$, where $f$ runs through all
$G_0$-invariant functions. Thus by Proposition~\ref{a40.70}
$T^{1,0}\O\subset(T^{1,0}\O)^{\bot\omega}$.  q.e.d.

\section{Main theorem: completeness of reductions for
generic coadjoint orbits}

We retain notations and conventions from two preceding
sections.

\abz\label{main}
\begin{th}
Let $U$ be an open set in  a $G$-orbit $M\subset\D$
such that $U'=U/G_0$ is a smooth manifold and
let $p:  U\rightarrow U'$ be the canonical projection.
Write $J$ for the holomorphic symplectic bihamiltonian
structure on $M$ associated with the restriction of the
standard holomorphic symplectic form $\omega=(c|_M)^{-1}$.

Then  the reduction $J'$ (cf. Corollary~\ref{a30.60} and
Definition~\ref{a30.70}) of $J|_U$ via $p$ is a complete
bihamiltonian structure at any $z'\not\in p(U\bigcap{\cal
R})$, where ${\cal R}$ is the Kronecker irregularity set
(see Definition~\ref{d40}).  In particular, $J'$ is
complete on $U'$ if $U$ is not contained in ${\cal R}$.

The realization $J$ of $J'$ is minimal (see Definition
\ref{a30.110}).
\end{th}

{\sc Proof.} We are going
to use Theorem \ref{ns} and the notations from it.

The foliation ${\cal K}$ of leaves of $p$ is a generic $CR$-foliation due to
Proposition \ref{a40.50}.  The
numbers $k_\la^z,k^z$ are constant in $z$ if $z'=p(z)$ is
fixed due to the $G_0$-invariance of all ingredients. Thus
the assumptions of Theorem~\ref{ns}
are satisfied.

By Corollaries \ref{CR} and \ref{a40.80} the number
$k^{z'}=\dim T^{1,0}_z{\cal K}\bigcap\linebreak
(T^{1,0}_z{\cal K})^{\bot\omega (z)}, z\in p^{-1}(z')$,
equals $r-\mu$, where $r=\ra\A$ and $\mu$ is the trivial
Kronecker dimension (see Definition~\ref{d40}), for any
$z'\in U'$.  We now shall prove that the number
$k^{z'}_\la=\dim T^\C{\cal K}\bigcap
(T^\C{\cal K})^{\perp\omega_\la (z)},z\in p^{-1}(z')$
satisfies the inequality
\sabz\label{f40.proof}
\begin{equation}
k^{z'}_\la\geq k^{z'}
\end{equation}
for
any $z'\in M',\la=(\la_1,\la_2)\in\C^2\smallsetminus\Lambda$.

For that purpose we shall put
$c_\la=\la_1\Re c+\la_2\Im c,\omega_\la=(c\la|_M)^{-1}$ and
use the fact that $(T^\C{\cal K})^{\bot\omega_\la}$, is
generated by the vector fields $c_\la(f)$, where $f$
 varies through the functions constant along ${\cal K}$.
Given a point $z\in U\smallsetminus{\cal R}$, we shall define
$r$ functions $g_{\la,1},\ldots,g_{\la,r}$ in a
neighbourhood $U_z$ of $z$ such that the vector fields
$c_\la(g_{\la,1}),\ldots,c_\la(g_{\la,r})$ are tangent to
${\cal K}$ and $r-\mu$ of them are independent on $U_z$.

Let us choose $U_z$ to be strongly admissible (see
Definition~\ref{d40.46} and Proposition~\ref{d40.55}) and
put
$$
g_{\la,1}=g_1(\tilde{\la}_2z+\tilde{\la}_1\bar{z}),\ldots
,g_{\la,r}=g_r(\tilde{\la}_2z+\tilde{\la}_1\bar{z}), $$
where
$(\tilde{\la}_1,\tilde{\la}_2)=(\frac{1}{2}(\la_1-\i\la_2),
\frac{1}{2}(\la_1+\i\la_2))\in \C^2$ is such that
$c_\la=\tilde{\la}_1c+\tilde{\la}_2\bar{c}$ and $g_1,\ldots
,g_r$ are functionally independent Casimir functions for
$c$.  These functions are functionally independent over
$U_z$ and their $\d$-differentials generate $(\ker
(\tilde{\la}_2c+\tilde{\la}_1\tilde{c})^{1,0}$ (see
Proposition~\ref{d40.50}). On the other hand, $\ker
c_\la=\ker c$ and Proposition~\ref{d40.60} implies that
over $U_z$ there are exactly $r-\mu$ independent among the
vector fields $c_\la(g_{\la,1}),\ldots,c_\la(g_{\la,r})$.

The nondegeneracy of $c_\la$ implies the independence of the vector fields
$c_\la(g_{\la,1}),\ldots,c_\la(g_{\la,r})$ at $z\in M\smallsetminus{\cal C}$. The
following equalities show that these vector fields are tangent to ${\cal K}$
($T{\cal K}$ is spanned by
$c(z_i)+\overline{c(z_i)},i=1,\ldots,n$,see
Proposition~\ref{d50})
\begin{eqnarray*}
c_\la(g_{\la,m})=\tilde{\la}_1c(g_{\la,m})+\tilde{\la}_2\bar{c}(g_{\la,m})=\\
\tilde{\la}_1\frac{\d g_{\la,m}}{\d z_i}c(z_i)+\tilde{\la}_2\frac{\d
g_{\la,m}}{\d \bar{z}_i}\bar{c}(\bar{z}_i)=\tilde{\la}_1\tilde{\la}_2\frac{\d
g_m}{\d z_i}|_{\tilde{\la}_2z+\tilde{\la}_1\bar{z}}(c(z_i)+\overline{c(z_i)}).
\end{eqnarray*}
Here we used the obvious identities
$$
\frac{\d g_{\la,m}}{\d
z_i}=\tilde{\la}_2\frac{\d g_m}{\d
z_i}|_{\tilde{\la}_2z+\tilde{\la}_1\bar{z}},$$

$$\frac{\d g_{\la,m}}{\d
\bar{z}_i}=\tilde{\la}_1\frac{\d g_m}{\d
z_i}|_{\tilde{\la}_2z+\tilde{\la}_1\bar{z}}.$$

Thus we have proved (\ref{f40.proof}) that is equivalent in
view of Theorem~\ref{ns} to the
following
$$
\ra c'_\la(z')\leq\ra c'(z'),z'\in U'\smallsetminus p({\cal R}).
$$
By the lower semi-continuity  of the function
$f(\la)= \ra c'_\la,\la\in\C^2$, this gives
$$
\ra c'_\la(z')=\ra c'(z'), z'\in U'\smallsetminus p({\cal R}).
$$
Thus we have obtained the constancy of $k^{z'}_\la$ in
$\la$ and the equality $k^{z'}_\la=k^{z'}$ for $z'\in
U'\smallsetminus p({\cal R})$. Since this number is also
independent of $z'$, condition $(iii)$ of Theorem \ref{ns}
is satisfied.

The minimality of the realization $J$ for $J'$ follows from
Corollary \ref{a40.80}. q.e.d.

We finish the paper by a characterization of the first
integrals (see Definition~\ref{20.05}) of the reduction
$J'$. It turns out that they are intimately related with
the first integrals of the canonical bihamiltonian
structure $\tilde{J}$ (see Section 5).

\abz\label{fin}
\begin{prop}
\begin{description}
\item[($i$)] Let $U\subset\D$ be a strongly admissible set
(see Definition~\ref{d40.46}) and let $F_1(U)$ be as in
Proposition~\ref{d40.57}. Then the functions from $F_1(U)$
are $G_0$-invariant.
\item[($ii$)] Retaining the hypotheses of
Theorem~\ref{main} assume that $U$ is strongly admissible.
Then the set of first integrals $F'_0(U')$ of the
bihamiltonian structure $J'$ over $U'$ is equal to the set
$(F_1(U))'$ consisting of the elements of the set $F_1(U)$
regarded as functions on $U'=U/G_0$.
\end{description}
\end{prop}

\noindent{\sc Proof.} ($i$) The elements of $Z_c^{hol}(U)$
are $G_0$-invariant by the definition, the $G_0$-invariance
of the functions from $\{\tilde{g};g\in Z_c^{hol}(U)\}$ is
established in the proof of Proposition~\ref{a40.80}. To
prove it for the elements from $\{g\circ\phi_\la|_U;g\in
Z_c^{hol}(U_\la)\}$, where $U_\la$ is an admissible
neighbourhood of $\phi_\la(U)$, it is sufficient to notice
the $G_0$-equivariance of the map $\phi_\la$.

\noindent ($ii$) It
follows from the proof of Theorem \ref{main} that the mentioned
there functions $g_{\la,1},\ldots,g_{\la,r}\in{\cal
F}_1(U)$ (now one can put $U_z=U$) being reduced generate
the space of the Casimir functions $Z_{c_\la'}$, where
$c_\la'$ is the reduction of the bivector $c_\la, \la\in
\C^2\smallsetminus\Lambda$.  Moreover, by the proof of
Proposition~\ref{a40.80} the functions
$\tilde{g_1},\ldots,\tilde{g_r}$ generate the space
$Z_{c'}$ after the reduction. q.e.d.

\abz\label{final}
\begin{cor}
Let $U\subset\D$ be a strongly admissible set. Then the
common level sets of functions from the family $\{\Re f,\Im
f;f\in F_1(U)\}$ form a foliation on  $U$ that is a
$CR$-lagrangian foliation (see Definition~\ref{a10.160}).
\end{cor}

\noindent{\sc Proof.} this foliation is the $CR$-lagrangian
foliation associated with the minimal realization
$(U,J|_U,{\cal K})$ of the complete bihamiltonian structure
$J'$ (see~\ref{a30.120}, \ref{a30.130}). q.e.d.

\end{document}